\documentclass{article}

\usepackage{arxiv}
\usepackage{amsmath}
\usepackage[utf8]{inputenc} % allow utf-8 input
\usepackage[T1]{fontenc}    % use 8-bit T1 fonts
\usepackage{hyperref}       % hyperlinks
\usepackage{url}            % simple URL typesetting
\usepackage{booktabs}       % professional-quality tables
\usepackage{amsfonts}       % blackboard math symbols
\usepackage{nicefrac}       % compact symbols for 1/2, etc.
\usepackage{microtype}      % microtypography
\usepackage{cleveref}       % smart cross-referencing
\usepackage{lipsum}         % Can be removed after putting your text content
\usepackage{graphicx}
\usepackage{doi}
\usepackage{algorithm}
\usepackage{algcompatible}
\usepackage{bm}
\usepackage{pifont}
\usepackage{amsthm}
\usepackage{verbatim}
\usepackage{amsopn}
\usepackage{physics}
\usepackage{tikz}
\usepackage{mathdots}
\usepackage{yhmath}
\usepackage{cancel}
\usepackage{color}
\usepackage{array}
\usepackage{multirow}
\usepackage{amssymb}
\usepackage{gensymb}
\usepackage{tabularx}
\usepackage{extarrows}
\usepackage{booktabs}
\usetikzlibrary{fadings}
\usetikzlibrary{patterns}
\usetikzlibrary{shadows.blur}
\usetikzlibrary{shapes}

\newcommand{\bx}{{\boldsymbol{x}}}
\title{Energy Dissipation Preserving Physics Informed Neural Network for  Allen-Cahn Equations}

\newif\ifuniqueAffiliation
\uniqueAffiliationtrue

\ifuniqueAffiliation % Standard variant of author block
\author{ \href{https://orcid.org/0000-0002-4546-519X}{\includegraphics[scale=0.06]{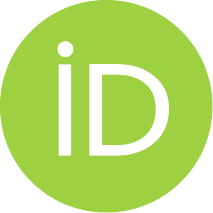}\hspace{1mm}Mustafa ~Kütük} \\
	Institute of Applied Mathematics\\
	Middle East Technical University\\
	06800, Ankara, Türkiye \\
	\texttt{mkutuk@metu.edu.tr} \\
	\And
	\href{https://orcid.org/0000-0002-0313-9767}{\includegraphics[scale=0.06]{orcid.pdf}\hspace{1mm}Hamdullah ~Yücel} \\
	Institute of Applied Mathematics\\
	Middle East Technical University\\
	06800, Ankara, Türkiye \\
	\texttt{yucelh@metu.edu.tr} \\
}
\else
\usepackage{authblk}

\newbox{\orcid}\sbox{\orcid}{\includegraphics[scale=0.06]{orcid.pdf}} 
\author[1]{%
	\href{https://orcid.org/0000-0002-4546-519X}{\usebox{\orcid}\hspace{1mm}Mustafa ~Kütük\thanks{\texttt{mkutuk@metu.edu.tr}}}%
}
\author[1]{%
	\href{https://orcid.org/0000-0002-0313-9767}{\usebox{\orcid}\hspace{1mm}Hamdullah ~Yücel\thanks{\texttt{yucelh@metu.edu.tr}}}%
}
\affil[1]{Institute of Applied Mathematics, Middle East Technical University, 06800, Ankara, Türkiye}
\fi

\renewcommand{\shorttitle}{\textit{arXiv} Template}

\hypersetup{
pdftitle={Energy Dissipation Preserving Physics Informed Neural Network for  Allen-Cahn Equations},
pdfsubject={q-bio.NC, q-bio.QM},
pdfauthor={Mustafa ~Kütük, Hamdullah ~Yücel},
pdfkeywords={First keyword, Second keyword, More},
}

\begin{document}
\maketitle

\begin{abstract}
    This paper investigates a numerical solution of Allen–Cahn equation with constant and degenerate mobility,  with polynomial and logarithmic energy functionals, with deterministic and random initial functions, and with advective term in  one, two, and three spatial dimensions, based on the physics-informed neural network (PINN). To improve the learning capacity of the PINN,  we incorporate the energy dissipation property of the Allen-Cahn equation  as a penalty term into the loss function of the network. To facilitate the learning process of random initials, we employ a continuous analogue of the initial random condition by utilizing the Fourier series expansion. Adaptive methods from traditional numerical analysis are also integrated  to enhance the effectiveness of the proposed PINN. Numerical results indicate a consistent decrease in the discrete energy, while also revealing phenomena such as phase separation and metastability.
\end{abstract}

% keywords can be removed

%%%%%%%%%%%%%%%%%%%%%%%%%%%%%%%%%%%%%%%%%%%%%%%%%%%%%%%%%%%%%%%
%%%%%%%%%%%%%%%%%%%%%%%%%%%%%%%%%%%%%%%%%%%%%%%%%%%%%%%%%%%%%%%%%%%%%%%%%%%%%%%%%%%%%%%%%%%%%%%%%%%%%%%%%%%%%%%%%%%%
%%%%%%%%%%%%%%%%%%%%%%%%%%%%%%%%%%%%%%%%%%%%%%%%%%%%%%%%%%%%%%%%%%%%%%%%%%%%%%%%%%%%%%%%%%%%%%%%%%%%%%%%%%%%%%%%%%%%
\section{Introduction}\label{sec:intro}

Phase field models  have a central role to understand the behaviour of the many complicated moving interface problems in material science \cite{ISteinbach_2009}, fluid dynamics \cite{JKim_2012}, fracture mechanics \cite{AKarma_DAKessler_HLevine_2001,BBourdin_GAFrancfort_JJMarigo_2000}, image analysis \cite{YLi_JKim_2011}, and mean curvature flow \cite{MBrassel_EBretin_2011}. Among these models, Allen-Cahn introduced in \cite{SMAllen_JWCahn_1979a} to describe the motion of anti-phase boundaries in crystalline solids at a fixed temperature is a particular case of gradient flows in the form of 
\begin{equation} \label{acgradient}
u_t = -\mu(u)\ \frac{\delta \mathcal{E}(u)}{\delta u},
\end{equation}
where $\frac{\delta \mathcal{E}(u)}{\delta u}$  represents the variational derivative of the free energy taken in the $L^2(\Omega)$-norm with  $\Omega \subset \mathbb{R}^d$ ($d=1,2,3$) as follows
\begin{equation}\label{energy}
\mathcal{E}(u)=\int_{\Omega} \left( \frac{\epsilon^2}{2}|\nabla u|^2 + F(u)  \right)d\bx.
\end{equation}
Specifically, the Allen-Cahn equation corresponds to a  nonlinear partial differential equation (PDE) in the form of 
\begin{equation}\label{allencahn}
u_t=\mu(u)(\epsilon^2\Delta u -f(u)), \quad (\bx,t) \in \Omega  \times (0,T],
\end{equation}
with the initial condition
\[
u(\bx,t=0) =u_0(\bx) \qquad \bx \in \Omega,
\]
and appropriate boundary conditions such as the Dirichlet, Neumann, or periodic boundary condition. Here, $u$  denotes the concentration of one the species of the alloy, known as the phase state between materials, on a bounded domain $\Omega \subset \mathbb{R}^d$ ($d=1,2,3$) with  the Lipschitz boundary $\partial \Omega$, $\epsilon$  represents the small interfacial length during the phase separation process, and $\mu(u)$ is the non--negative mobility function. 
Moreover, the nonlinear $f(u)$ is the derivative of a free energy functional $F(u)$, which is characterized by logarithmic or polynomial functions. The logarithmic free energy potential (see, e.g., \cite{JWBarrett_JFBlowey_HGarcke_2000}) can be formulated as 
\begin{equation} \label{Flory_Huggins_Energy}
F(u) = \frac{\theta}{2} [(1+u) \ln(1+u) + (1-u) \ln(1-u)] - \frac{\theta_c}{2} u^2, \,\,\,\,\, 0<\theta\leq\theta_c, 
\end{equation}   
where $\theta$ and $\theta_c$ are absolute and transition temperatures, respectively. The derivative of the logarithmic free energy potential is then equivalent to 
\begin{equation*} 
f(u) = F'(u) = \frac{\theta}{2} \ln\Bigg[\frac{1+u}{1-u}\Bigg] - \theta_c u. 
\end{equation*}  
On the other hand, the polynomial free energy potential, also called as quartic double--well, is an approximation of the logarithmic ones when the temperature $\theta$ closes to $\theta_c$, formulated as follows
\begin{equation} \label{Polynomial_Energy}
F(u) = \frac{1}{4}(1 - u^2)^2,
\end{equation} 
and its derivative is $f(u) = u^3-u$. It is well-known that the solution $u(\bx,t)$ of the Allen–Cahn equation has two crucial properties; the decay of the total free energy  for $\mu(u)>0$, that is,
\begin{equation}\label{energy_law}
 \frac{d \mathcal{E}(u(t))}{dt} \leq 0,
\end{equation}
which is a typical property of the gradient flows, and the maximum bound principle. In the double-well potential \eqref{Polynomial_Energy}, the solution always stays inside the interval $[-1, 1]$, while the solution of the logarithmic potential \eqref{Flory_Huggins_Energy} is bounded in the interval $[-s,s]$  for all the time, where $1>s>0$ is the constant satisfying $f(s)=0$.

The Allen-Cahn equation \eqref{allencahn} has been the subject of a great deal of research, primarily because of its relationship with intriguing and complex geometric problems involving moving surfaces; see, e.g. \cite{MJWard_1996,PdeMottoni_MSchatzma_1991,LBronsard_RVKohn_1991} and references therein. As $\epsilon$ tends towards zero, the solution of the Allen-Cahn equation \eqref{allencahn} behaves like a piecewise constant with values  $\pm 1$, in much of two bulk regions which are separated by a diffusive interfacial layer of thickness $\mathcal{O}(\epsilon)$. For finite but small $\epsilon$, the solution within this interfacial layer remains smooth but it develops a large spatial gradient. Such layers move in time, mimicking dynamically evolving fronts. Because of that, the Allen-Cahn equation has been considered as one of the fundamental model equations for the diffuse interface approach developed to study phase transitions and interfacial dynamics in materials science \cite{LQChen_2002}. In \cite{SMAllen_JWCahn_1979a}, it was first formally proved that the zero level set of the solution $u$ converges  to a surface denoted by $\Gamma_t$ which evolves according to the geometric law $V = \kappa$, where $V$ represents the normal velocity of the surface $\Gamma_t$ and $\kappa$ is its mean curvature. Later, Evans et al. in \cite{LCEvans_HMSoner_PESouganidis_1992} rigorously justified this limit by proving a comprehensive result: for $t \geq 0$, the limit of the zero level set of the solution $u$ is encompassed within the generalized solution of the motion by mean curvature flow, as established in  \cite{YGChen_YGiga_SGoto_1989,LCEvans_JSpruck_1991}. Then, it was also shown in \cite{TIlmanen_1993} that this limit is one of Brakke's motion by mean curvature solution.

Clearly, the Allen-Cahn equation \eqref{allencahn} is crucial for understanding phase transitions in materials science and examining curvature-driven flows in geometry. It is especially valuable when singularities arise in these flows, as the solution $u$ persists over time, requiring minimal adjustments to manage potential singularities. While analytical methods are used to understand the dynamics of the dynamics of the solutions, explicit analytic solutions are generally not possible due to the nonlinearity of the equation. Therefore, it becomes essential to find an efficient and accurate numerical solution. However, it is not an easy task due to 
the existence of a nonlinear term and  of the small interfacial length parameter $\epsilon$. An extensive numerical investigation of the Allen-Cahn equations has been studied  by using different numerical techniques in the spatial domain, such as finite element methods in \cite{FLiu_JShen_2013a},  discontinuous Galerkin methods in \cite{XFeng_YLi_2015a,BKarasozen_MUzunca_ASFilibelioglu_Hyucel_2018}, and in the temporal domain,  such as
the implicit-explicit (IMEX) techniques \cite{JShen_XYang_2010a,XFeng_HSong_TTang_JYang_2013a}, the average vector field (AVF) method
\cite{BKarasozen_MUzunca_ASFilibelioglu_Hyucel_2018,ECelledoni_VGrimm_RIMclachlan_DIMclaren_DOneale_BOwren_GRWQuispel_2012a}, and the splitting methods 
\cite{JShen_CWang_XWang_SWise_2012}. In addition to the discretization approaches, analytical approximate techniques  have been used to solve such kind of nonlinear PDEs, such as Adomian decomposition method, homotopy perturbation method, homotopy analysis method, tanh-function method, (variational iteration) method; see, e.g., \cite{TKAkinfe_ACLoyinmi_2022,SHussain_FHaq_AShah_DAbduvalieva_AShokri_2024}. In these studies, the primary goal is to find a numerical scheme that satisfies the energy dissipative property \eqref{energy_law} and adheres to the maximum bound principle, in addition to obtaining (optimal) a priori error estimates.

%JDEyre,

%As $\epsilon$ tends towards zero, it forms boundaries between the two stable \textcolor[rgb]{0.00,0.50,1.00}{phases} $u=\pm 1$, which is a process known as phase separation.

The aforementioned numerical methods can be accurate up to the given threshold but they have some restrictions such as mesh dependence, high computational burden for the nonlinear PDEs. On the other hand, in recent years, the usage of deep neural networks (DNNs)  has lead to a significant achievement in the several areas, such as visual recognition \cite{AKrizhevsky_ISutskever_GHinton_2012}, cognitive science \cite{BMLake_RSalakhutdinov_JBTenenbaum_2015} as well as solving differential equations  \cite{JBerg_KNystrom_2018,WE_BYu_2018,BLi_STang_HYu_2020,LLu_XMeng_ZMao_GEKarniadakis_2021,MRaissi_PPerdikaris_GEKarniadakis_2019,JSirignano_KSpiliopoulos_2018,YLiao_PMing_2021}. Among them, Physics-Informed Neural Network (PINN) introduced by Raissi et al. in \cite{MRaissi_PPerdikaris_GEKarniadakis_2019} has received  great attention thanks to its flexibility in tackling a wide range of forward and inverse problems involving PDEs. The idea behind PINN is based on the universal approximation theorem \cite{GCybenko_1989,KHornik_MStinchcombe_HWhite_1989}, which asserts that a neural network can approximate any function given sufficient complexity. In this structure,  weights and biases of the neural network model are  optimized according to the loss function, containing the physics of the underlying problem, which are basically governing equations, boundary conditions, and initial condition. Up to now, the PINN has been used to solve  different types of problems in computational science and engineering, such as inverse problems \cite{YChen_LLu_GEKarniadakis_LDNegro_2020}, fluid dynamics \cite{ZMao_ADJagtap_GEKarniadakis_2020,MSikora_PKrukowski_APaszynska_MPaszynski_2024}, parameter estimation \cite{MDaneker_ZZhang_GEKarniadakis_LLu_2023}, topology optimization \cite{LLu_RPestourie_WYao_ZWang_FVerdugo_SGJohnson_2021}, fractional PDEs \cite{GPang_LLu_GEKarniadakis_2019}, and stochastic PDEs \cite{DZhang_LGuo_GEKarniadakis_2018,SKarumuri_RTripathy_IBilionis_JPanchal_2020}. In addition, different variants of the standard PINN have been proposed to increase the performance of PINN, such as meta-learning \cite{AFPsaros_KKawaguchi_GEKarniadakis_2022}, gradient-enhanced \cite{JYu_LLu_XMeng_GEKarniadakis_2022}, balance of weights \cite{SWang_YTeng_PPerdikaris_2021}, decomposition of spatial domain \cite{ADJagtap_GEKarniadakis_2020}, adaptive sampling \cite{CWu_MZhu_QTan_YKartha_LLu_2023}, adaptivity in temporal domain  \cite{CLWight_JZhao_2021,SWang_SSankaran_PPerdikaris_2024}, adaptive activation function \cite{ADJagtap_KKawaguchi_GEKarniadakis_2020}, and enforcing boundary conditions \cite{SDong_NNi_2021}. Although the PINN has been widely employed for approximating PDEs, theoretical investigations into their convergence and error analysis remain limited. In \cite{YShin_JDarbon_GEmKarniadakis_2020}, the authors explain why the PINN works and shows its consistency for the linear elliptic and parabolic PDEs. An estimation of the generalization error by means of the training error and the number of training data points is provided in \cite{SMishra_RMolinaro_2021,SMishra_RMolinaro_2022} for forward PDEs and inverse problems. Total error is  bounded in terms of training samples, depth and width of the deep neural networks in various studies:  \cite{TDeRyck_ADJagtap_SMishra_2023} for Navier-Stokes equation,  \cite{YJiao_YLai_DLi_XLu_FWang_YWang_JZYang_2022} for the linear elliptic and parabolic PDEs,  \cite{GZhang_JLin_QZhai_HYang_XChen_XZheng_ITakLeong_2024} for phase field problems, and  \cite{YQian_YZhang_YHuang_SDong_2023} for hyperbolic PDEs. The mathematical foundation for PINNs in approximating PDE solutions continues to be an active area of research. However, a theoretical discussion is not the scope of the present study and  can be considered as a future study. Instead, we mainly focus on  presenting a novel PINN framework to solve a class of Allen-Cahn equations. Compared to the conventional methods, such as finite difference, finite element, deep learning approaches like PINN is mesh-free thanks to the automatic differentiation and can avoid the curse of dimensionality; see, e.g., \cite{TPoggio_HMhaskar_LRosasco_BMiranda_QLiao_2017}. However, for low  dimensional PDEs, it is still not straightforward to claim that the computational accuracy obtained from DNNs is better than the ones obtained by the conventional methods; see, e.g., \cite{TGGrossmann_UJKomorowska_JLatz_CBSchonlieb_2024} for more discussion.

Although the standard PINN (std-PINN) has been widely accepted and has yielded remarkable results across a range of problems in the computational science and engineering, it is not always capable of solving the Allen-Cahn equation, a typical example of phase field models, due to the sharp transition layers, evolution of the layers in time, high sensitivity to the initial conditions, or the choice of hyperparameters, including the learning rate, activation function, optimizer, loss function weights, and the width and depth of the neural network. To enhance the efficiency of PINN in solving Allen-Cahn type PDEs, the authors in \cite{YGeng_YTeng_ZWang_LJu_2024} use two specially designed convolutional neural networks (CNNs) and the loss functions correspond to the full-discrete systems obtained from finite difference methods in both space and time. Similarly, the multi-step discrete time models  with adaptive collocation strategy are considered in  \cite{HXu_JChen_FMa_2022}. The authors in  \cite{CLWight_JZhao_2021} involve the idea of adaptivity in both time and space by sampling the data points, while in \cite{RMattey_SGhosh_2022}, the same neural network is retrained  over successive time segments, while satisfying the obtained solution for all previous time segments. In addition, the system of the  Allen-Cahn equation \eqref{allencahn} first reduced into a first-order problem and then the converted minimization problem is
approximated by using a deep learning approach in \cite{ASingh_RKSinha_2023}.  A theoretical perspective for the propagation of  PINN errors is also given in \cite{GZhang_JLin_QZhai_HYang_XChen_XZheng_ITakLeong_2024} for the Allen-Cahn  equations. A recent study by Guo et al. in \cite{JGuo_HWang_CHou_2024}  introduces an adaptive energy-based sequential training method for the PINN approach to enhance their performance in solving the Allen-Cahn equation, particularly for long-duration simulations. In the aforementioned works, the proposed PINN structures are usually tested on the Allen-Cahn equations with polynomial potentials, while to the best of our knowledge there is no PINN work on the Allen-Cahn equations with logarithmic potentials, where the nonlinearity plays a more critical role. In addition, the behavior of the energy functional, which is physically more crucial for such kind of problems, is often neglected. Contrary to the previous studies, we here propose a novel methodology based on preserving of the energy dissipation to predict the dynamics of the Allen-Cahn equation with constant and degenerate mobility,  with polynomial and logarithmic energy functionals, with deterministic and random initial functions, and with advective term in  one, two, and three spatial dimensions. The proposed network approach guarantees the decay of energy and also plays a key role in accurately learning the dynamics of the Allen-Cahn equation. Embedding of different conservation constraints, such as mass and momentum conservation, into the PINN architecture are also considered in \cite{GZWu_YFang_NAKudryashov_YYWang_CQDai_2022,SLin_YChen_2022} for different types of PDEs.

Our specific contributions can be summarized as:
\begin{itemize}
    \item We propose a novel PINN approach  based on preserving of the energy dissipation to learn the dynamics of the Allen-Cahn equation more accurately.
    \item We offer a detailed set of benchmark examples to evaluate the performance of the proposed approach. These examples include logarithmic and polynomial free energy potentials, constant and degenerate mobility functions, deterministic and random initial functions,  and advective term in one, two, and three spatial dimensions.    
\end{itemize}

The rest of this manuscript is outlined as follows: In next section, we briefly review the standard physics-informed neural network (std-PINN). In Section~\ref{sec:energy_pinn}, we introduce our main contribution, which is the addition of energy dissipation constraint  \eqref{energy_law}  into the loss function. Section~\ref{sec:algorithm} presents numerical strategies utilizing adaptive approaches  to enhance the performance of the PINN. Numerical simulations of  various benchmark problems are provided in Section~\ref{sec:numeric}. Last, we give some concluding remarks in Section~\ref{sec:conclusion} based on the findings in the paper.

%%%%%%%%%%%%%%%%%%%%%%%%%%%%%%%%%%%%%%%%%%%%%%%%%%%%%%%%%%%%%%%%%%%%%%%%%%%%%%%%%%%%%%%%%%%%%%%%%%%%%%%%%%%
\section{Physics-informed neural network}\label{sec:pinn} 

In this section, we first review the nonlinear mapping in the deep learning architecture \cite{IGoodfellow_YBengio_ACourville_2017},  and then the standard physics-informed neural network (std-PINN) introduced in \cite{MRaissi_PPerdikaris_GEKarniadakis_2019}.

Let  a neural network of depth $L$, that is, $\mathcal{N}^L(\boldsymbol{x}): \, \mathbb{R}^{N_0} \rightarrow \mathbb{R}^{N_L}$ consist of one input layer,  $L-1$ hidden layer, and one output layer. Each layer which contains $N_{\ell}$ neurons transmits the output $\boldsymbol{x}^{\ell}$ to the $(\ell +1)$th layer as the input data. The connection between layers is constructed by an affine transformation $\mathcal{T}$ and an activation function $\sigma(\cdot)$
\begin{equation*}
\boldsymbol{x}^{\ell} = \sigma \big( \mathcal{T}_{\ell} (\boldsymbol{x}^{\ell -1}) \big)  = \sigma \big( w^{\ell - 1} \boldsymbol{x}^{\ell - 1} + b^{\ell - 1}\big),
\end{equation*}
where $w^{\ell - 1} \in \mathbb{R}^{N_{\ell} \times N_{\ell -1}}$  and $b^{\ell - 1} \in \mathbb{R}^{N_{\ell - 1}}$ denote the weight and bias of the $\ell$th layer, respectively. In general, linear functions are often inadequate for capturing the complexity of a problem. Therefore, nonlinear functions such as ReLU, sigmoid, $\tanh$, swish, etc. \cite{TSzandala_2021} are commonly used as activation functions. Hence, for a given  input $\boldsymbol{x}^0$, we have
\begin{equation*}
\widehat{u}(\boldsymbol{x}^{0}, \theta) = \big( \sigma \circ \mathcal{T}_{L} \circ \sigma \circ \mathcal{T}_{L-1} \circ  \cdots \circ \sigma \circ   \mathcal{T}_{1} \big)(\boldsymbol{x}^0),
\end{equation*}
where  $\widehat{u}(\cdot,\cdot )$ is the output of the learning process and  $\theta= \{ w^{\ell}, b^{\ell}  \}_{\ell=1}^{L}$ represent the trainable parameters. Typically, the same activation function is applied to all hidden layers. However, the activation function of the output layer can differ from those of the hidden layers, depending on the structure of the problem. Before training a neural network, $\theta$ is need to be initialized by  for instance, Xavier initialization \cite{XGlorat_Y.Bengio_2010}, He initialization \cite{KMHe_XYZhang_SQRen_JSun_2015}, etc.

Let $\Omega \subset \mathbb{R}^d$ ($d=1,2,3)$ be an open bounded domain in $\mathbb{R}^d$  with a sufficiently smooth boundary $\partial \Omega$, and $(0,T]$ be a time interval with $T< + \infty$. Then, a general form of a $m$th order partial differential equation is given in the form of
\begin{eqnarray}\label{PDE}
u_t  + \mathcal{N} (u, u_{\bx}^{(1)},u_{\bx}^{(2)}  \ldots, u_{\bx}^{(m)})) &=& 0 \quad \bx \in \Omega \subset \mathbb{R}^d, \;  t \in (0,T],
\end{eqnarray}
where $\bx$ and $t$ denote the spatial component with $\Omega \subset \mathbb{R}^d$ and time, respectively. $\mathcal{N}(\cdot)$  is a nonlinear function of the solution $u(\bx,t)$ and its spatial derivatives $\big( u_{\bx}^{(1)}(\bx,t),u_{\bx}^{(2)}(\bx,t)  \ldots, u_{\bx}^{(m)}(\bx,t) \big)$. Corresponding boundary and initial conditions are given, respectively, by
\begin{eqnarray}\label{bcic}
\mathcal{B}(u, u_{\bx}^{(1)}, , \ldots, u_{\bx}^{(m-1)}) &=& G(\bx,t) \qquad (\bx,t) \in \partial \Omega \times (0,T], \label{bc} \\
u(\bx,0) &=& H(\bx) \quad \qquad \bx \in \Omega, \label{ic}
\end{eqnarray}
where the operator $\mathcal{B}(\cdot)$ allows Dirichlet, Neumann, periodic, and also higher order boundary conditions. Our aim is to develop a numerical solution for the PDE problem \eqref{PDE} with the boundary condition \eqref{bc} and the initial condition \eqref{ic} using the neural networks.

In the prediction process of std-PINN, the spatial and temporal points constitute the inputs of network and  randomly initialized parameters are updated by minimizing the loss function in the PINN formulation, including  three error components defined on the residual, initial condition, and boundary conditions in the sense of
the mean--squared error (MSE). Denoting $\widehat{u}(\bx,t)$ as the output of neural network, the components of the total loss function are introduced as follow:
\begin{itemize}
	\item MSE on the residual PDE
	\begin{equation}
	\hbox{loss}_r =  \frac{1}{n_r} \sum \limits_{k=1}^{n_r} \Big( \underbrace{ \widehat{u}_t (\bx_k^r,t_k^r)  + \mathcal{N} (\widehat{u}(\bx_k^r,t_k^r), \ldots, \widehat{u}_{\bx}^{(m)} (\bx_k^r,t_k^r))}_{R(\bx_k^r,t_k^r)} \Big)^2,
	\end{equation}
	where $(\bx_k^r, t_k^r) \in \Omega \times (0,T]$  and $n_r$ denote the interior collocation points and its size, respectively. 
	\item MSE on the initial condition
	\begin{equation}
	\hbox{loss}_i =  \frac{1}{n_i} \sum \limits_{k=1}^{n_i} \Big( \underbrace{\widehat{u}(\bx_k^i,0) - H(\bx_k^i)}_{I(\bx_k^i)} \Big)^2, 
	\end{equation}
	where  $n_i$ is the number of initial training data $x_k^i \in \Omega$.
	\item MSE on the boundary conditions
	\begin{equation}
	\hbox{loss}_b =  \frac{1}{n_b} \sum \limits_{k=1}^{n_b} \Big( \underbrace{\mathcal{B}(\widehat{u}(\bx_k^b,t_k^b), \ldots, \widehat{u}_{x}^{(m-1)}(\bx_k^b,t_k^b)) -   G(\bx_k^b,t_k^b)}_{B(\bx_k^b,t_k^b)} \Big)^2,
	\end{equation}
	where  $(\bx_k^b,t_k^b) \in \partial \Omega  \times (0,T]$  represents the set of training points to compute MSE on the boundary for various times and  $n_b$ is the number of training points on the boundary.
\end{itemize}

In the PINN structure, the derivatives of the network with respect to time $t$ and space $\bx$ are computed by the automatic differentiation \cite{AGBaydin_BAPearlmutter_AARadul_JMSiskind_2018}. On the other hand, the choice of sampling strategy of collocation points $(\bx_k, t_k)$ in the training of PINNs is crucial for enhancing both the accuracy and computational efficiency. Typically, in many PINN studies, the collocation points are determined at the start of the training and remain unchanged throughout the process. In this context, equispaced uniform grids and uniformly random sampling are two basic sampling methods. However, there are also other techniques, such as Latin hypercube sampling (LHS) \cite{MDShields_JZhang_2016} and Sobol sampling \cite{MRenardy_LRJoslyn_JAMillar_DEKirschner_2021}. The LHS is a stratified Monte Carlo sampling technique that produces random samples distributed within intervals, ensuring equal probability and a normal distribution for each range. In contrast, the Sobol sequence, a type of quasirandom low-discrepancy sequence, is often utilized as an alternative to uniformly distributed random numbers and typically yield superior results in various applications, including numerical integration. Instead of utilizing fixed residual points throughout the training process, one could opt to select a new set of residual points at every specified optimization iteration; see, e.g., \cite{LLu_XMeng_ZMao_GEKarniadakis_2021}.   
Further, recent studies, such as  \cite{LLu_XMeng_ZMao_GEKarniadakis_2021,CLWight_JZhao_2021,JMHanna_JVAguadoa_SCCardonab_RAskri_DBorzacchielloa_2022} 
 have shown that adaptive sampling strategies significantly enhance the distribution of residual points during training, leading to improved accuracy; see Section~\ref{sec:algorithm} for a detailed discussion.

Hence, the total loss function of the std-PINN is given by adding all the aforementioned mean--square errors
\begin{equation}\label{lossEq}
\hbox{Loss}(\theta) = \lambda_r \, \hbox{loss}_r(\theta) + \lambda_i \, \hbox{loss}_i(\theta) +  \lambda_b \, \hbox{loss}_b(\theta),
\end{equation}
where $\theta$ is the set of trainable parameters, consisting of the weights and biases, and the parameters $\lambda_r$, $\lambda_i$, and $\lambda_b$ scale the loss terms. Then, one wishes to find the optimal value of training parameters $\theta$  by minimizing the loss function  (\ref{lossEq}) until a given threshold tolerance value or  a prescribed maximum number of iterations. In the literature, there are several effective  optimization techniques, such as stochastic gradient descent (SGD), L-BFGS, ADAM, etc., see, e.g., \cite{LBottou_FECurtis_JNocedal_2018a} and references therein for more details.  

%%%%%%%%%%%%%%%%%%%%%%%%%%%%%%%%%%%%%%%%%%%%%%%%%%%%%%%%%%%%%%%%%%%%%%%%%%%%%%%%%%%%%%%%%%%%%%%%
\section{Penalization of the energy dissipation}\label{sec:energy_pinn}

To achieve an accurate and efficient numerical approximation for the solutions of nonlinear PDEs, it is essential to maintain the qualitative structures of the underlying system, such as mass conservation and energy dissipation. Motivated by this, we propose a new methodology that focuses on preserving energy dissipation to predict the dynamics of the Allen-Cahn equation. 

%Now, if we consider a general optimization problem in the form of 
%\begin{equation*}
%   \min \limits_{\bx \in \mathbb{R}^n} \,\, f(\bx) \quad \hbox{subject to} \quad c(\bx) \leq 0,
%\end{equation*}
%where $c(\bx)$  is the inequality constraint of the optimization problem. From the optimization theory \cite{JNocedal_SJWright_2006a}, we can penalize it as %follows
%\begin{equation*}
%    \min \limits_{\bx \in \mathbb{R}^n} f(\bx ) + \frac{\sigma}{2} \|\max (0, c(\bx))\|_2^2,
%\end{equation*}
%where $\sigma$ is a penalty parameter. 

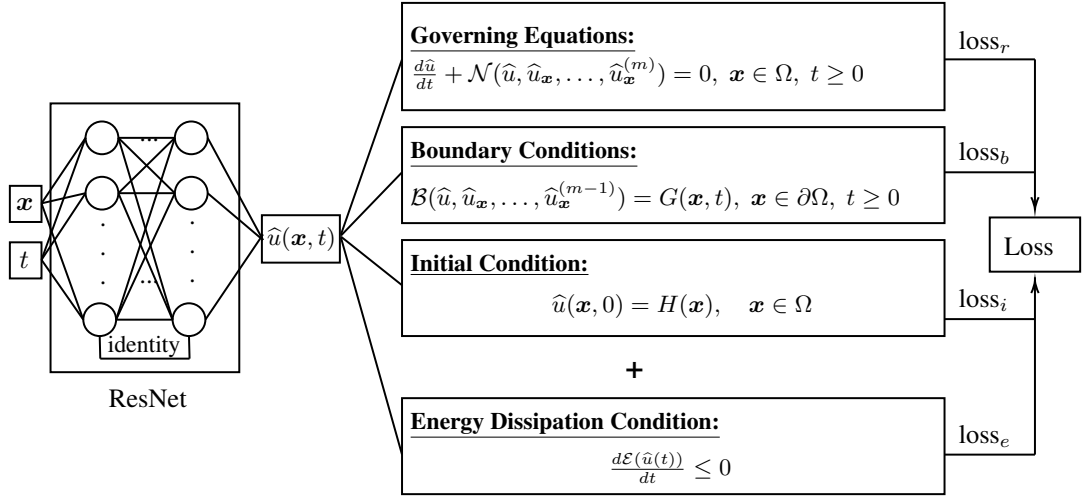
\begin{figure}[htp!]
    \begin{center}
\tikzset{every picture/.style={line width=0.75pt}} %set default line width to 0.75pt        

\begin{tikzpicture}[x=0.75pt,y=0.75pt,yscale=-1,xscale=1,scale=0.57]
%uncomment if require: \path (0,453); %set diagram left start at 0, and has height of 453

%Shape: Square [id:dp9023623002429422] 
\draw   (23,122) -- (51,122) -- (51,154) -- (23,154) -- cycle ;
%Shape: Square[id:dp9615699278695689] 
\draw   (23,172) -- (51,172) -- (51,204) -- (23,204) -- cycle ;
%Shape: Circle [id:dp6445064646553267] 
\draw   (90,79.5) .. controls (90,71.49) and (96.49,65) .. (104.5,65) .. controls (112.51,65) and (119,71.49) .. (119,79.5) .. controls (119,87.51) and (112.51,94) .. (104.5,94) .. controls (96.49,94) and (90,87.51) .. (90,79.5) -- cycle ;
%Shape: Circle [id:dp09696058407406216] 
\draw   (90,128.5) .. controls (90,120.49) and (96.49,114) .. (104.5,114) .. controls (112.51,114) and (119,120.49) .. (119,128.5) .. controls (119,136.51) and (112.51,143) .. (104.5,143) .. controls (96.49,143) and (90,136.51) .. (90,128.5) -- cycle ;
%Shape: Circle [id:dp9424574735616178] 
\draw   (88,240.5) .. controls (88,232.49) and (94.49,226) .. (102.5,226) .. controls (110.51,226) and (117,232.49) .. (117,240.5) .. controls (117,248.51) and (110.51,255) .. (102.5,255) .. controls (94.49,255) and (88,248.51) .. (88,240.5) -- cycle ;

%Straight Lines [id:da6571537312333078] 
\draw    (51,137.5) -- (90,79.5) ;
%Straight Lines [id:da19327733417528337] 
\draw    (51,137.5) -- (90,128.5) ;
%Straight Lines [id:da31095156706631166] 
\draw    (51,187.5) -- (88,240.5) ;
%Straight Lines [id:da958228725850468] 
\draw    (51,137.5) -- (88,240.5) ;
%Straight Lines [id:da17987799972769758] 
\draw    (51,187.5) -- (90,128.5) ;
%Straight Lines [id:da6018410427843461] 
\draw    (51,187.5) -- (90,79.5) ;
%Shape: Circle [id:dp8602582774428462] 
\draw   (169,79.5) .. controls (169,71.49) and (175.49,65) .. (183.5,65) .. controls (191.51,65) and (198,71.49) .. (198,79.5) .. controls (198,87.51) and (191.51,94) .. (183.5,94) .. controls (175.49,94) and (169,87.51) .. (169,79.5) -- cycle ;
%Shape: Circle [id:dp5175610341834116] 
\draw   (169,128.5) .. controls (169,120.49) and (175.49,114) .. (183.5,114) .. controls (191.51,114) and (198,120.49) .. (198,128.5) .. controls (198,136.51) and (191.51,143) .. (183.5,143) .. controls (175.49,143) and (169,136.51) .. (169,128.5) -- cycle ;
%Shape: Circle [id:dp45501067314168564] 
\draw   (167,240.5) .. controls (167,232.49) and (173.49,226) .. (181.5,226) .. controls (189.51,226) and (196,232.49) .. (196,240.5) .. controls (196,248.51) and (189.51,255) .. (181.5,255) .. controls (173.49,255) and (167,248.51) .. (167,240.5) -- cycle ;
%Straight Lines [id:da3888984524220154] 
\draw    (119,79.5) -- (169,79.5) ;
%Straight Lines [id:da7229707152468903] 
\draw    (119,79.5) -- (169,128.5) ;
%Straight Lines [id:da11237826943461826] 
\draw    (119,79.5) -- (167,240.5) ;
%Straight Lines [id:da2661226241487449] 
\draw    (119,128.5) -- (169,79.5) ;
%Straight Lines [id:da6953103744215643] 
\draw    (119,128.5) -- (169,128.5) ;
%Straight Lines [id:da30978954021826755] 
\draw    (119,128.5) -- (167,240.5) ;
%Straight Lines [id:da6351888937718684] 
\draw    (167,240.5) -- (117,240.5) ;
%Straight Lines [id:da008961702641864688] 
\draw    (117,240.5) -- (169,128.5) ;
%Straight Lines [id:da10300279744914365] 
\draw    (117,240.5) -- (169,79.5) ;
%Shape: Square [id:dp4382709365209154] 
\draw   (246,148) -- (315,148) -- (315,190) -- (246,190) -- cycle ;
%Straight Lines [id:da9362401310176978] 
\draw    (198,79.5)  -- (247,170) ;
%Straight Lines [id:da5674878003594039] 
\draw    (198,128.5) -- (247,170) ;
%Straight Lines [id:da7180026475465515] 
\draw    (196,240.5) -- (247,165) ;
%Shape: Rectangle [id:dp9619708158556133] 
\draw   (59,49) -- (226,49) -- (226,287) -- (59,287) -- cycle ;
%Straight Lines [id:da1434033909732575] 
\draw    (102.5,255) -- (103,275) ;
%Straight Lines [id:da9886716194668788] 
\draw    (103,275) -- (182,275) ;
%Straight Lines [id:da4638721638554635] 
\draw    (181.5,255) -- (182,275) ;

% Text Node
\draw (26,130) node [anchor=north west][inner sep=0.75pt]   [align=left] {$\bx$};
% Text Node
\draw (29,178) node [anchor=north west][inner sep=0.75pt]   [align=left] {$t$};
% Text Node
\draw (100,153) node [anchor=north west][inner sep=0.75pt]   [align=left] {.\\.\\.};
% Text Node
\draw (179,151) node [anchor=north west][inner sep=0.75pt]   [align=left] {.\\.\\.};
% Text Node
\draw (135,203) node [anchor=north west][inner sep=0.75pt]   [align=left] {...};
% Text Node
\draw (248,156) node [anchor=north west][inner sep=0.75pt]   [align=left] {\footnotesize $\displaystyle \widehat{u}(\bx,t)$};
% Text Node
\draw (107,252) node [anchor=north west][inner sep=0.75pt]   [align=left] {\footnotesize identity};
% Text Node
\draw (108,300) node [anchor=north west][inner sep=0.75pt]   [align=left] {ResNet};
% Text Node
\draw (135,77) node [anchor=north west][inner sep=0.75pt]   [align=left] {...};

%Shape: Rectangle [id:dp9072472866320778] 
\draw   (370,-40) -- (850,-40) -- (850,55) -- (370,55) -- cycle ;
% Text Node
\draw (375,-22) node [anchor=north west][inner sep=0.75pt]    
{\footnotesize \underline{\textbf{Governing Equations:}}};

\draw (375, 5) node [anchor=north west][inner sep=0.75pt]    
{\footnotesize $ \frac{d\widehat{u}}{d t}  + \mathcal{N} (\widehat{u}, \widehat{u}_{\bx}, \ldots, \widehat{u}_{\bx}^{(m)}) = 0, \; \bx \in \Omega, \;  t \geq 0$};

%Shape: Rectangle [id:dp9072472866320778] 
\draw   (370,70) -- (850,70) -- (850,155) -- (370,155) -- cycle ;

% Text Node
\draw (375,80) node [anchor=north west][inner sep=0.75pt]    
{\footnotesize \underline{\textbf{Boundary Conditions:}}};
\draw (375,115) node [anchor=north west][inner sep=0.75pt]    {\footnotesize $\mathcal{B}(\widehat{u}, \widehat{u}_{\bx},  \ldots, \widehat{u}_{\bx}^{(m-1)}) = G(\bx,t), \; \bx \in \partial \Omega, \;  t \geq 0$};

%Shape: Rectangle [id:dp9072472866320778] 
\draw   (370,170) -- (850,170) -- (850,255) -- (370,255) -- cycle ;

% Text Node
\draw (375,180) node [anchor=north west][inner sep=0.75pt]    
{\footnotesize\underline{\textbf{Initial Condition:}}};

\draw (500,215) node [anchor=north west][inner sep=0.75pt]    {\footnotesize $\widehat{u}(\bx,0) = H(\bx), \quad \bx \in \Omega$};

% Text Node
\draw (565,275) node [anchor=north west][inner sep=0.75pt]   [align=left] {\large \textbf{+}};

%Shape: Rectangle [id:dp7919141862986676] 
\draw   (370,310) -- (850,310) -- (850,395) -- (370,395) -- cycle ;
% Text Node
\draw (375,320) node [anchor=north west][inner sep=0.75pt]    
{\footnotesize \underline{\textbf{Energy Dissipation Condition:}}};

\draw (550,355) node [anchor=north west][inner sep=0.75pt]    {\footnotesize $\frac{d \mathcal{E}(\widehat{u}(t))}{dt} \leq 0$};

%Straight Lines [id:da8662181883435198] 
\draw    (316,165) -- (370,10) ;
%Straight Lines [id:da7108736886903582] 
\draw    (316,166) -- (370,110) ;
%Straight Lines [id:da582018451161687] 
\draw    (316,167) -- (370,210) ;
%Straight Lines [id:da8279704207397627] 
\draw    (316,168) -- (370,360) ;

%Straight Lines [id:da16121947632070066] 
\draw    (850,10) -- (930,10) ;
\draw (860,-18) node [anchor=north west][inner sep=0.75pt]    {$\hbox{loss}_r $};

%Straight Lines [id:da7849952212679991] 
\draw    (850,110) -- (930,110) ;
\draw (860,81) node [anchor=north west][inner sep=0.75pt]    {$\hbox{loss}_b $};

%Straight Lines [id:da10461000239110296] 
\draw    (850,240) -- (930,240) ;
\draw (860,211) node [anchor=north west][inner sep=0.75pt]    {$\hbox{loss}_i $};

%Straight Lines [id:da9789761056609332] 
\draw    (850,360) -- (930,360) ;
\draw (860,331) node [anchor=north west][inner sep=0.75pt]    {$\hbox{loss}_e $};

%Shape: Rectangle [id:dp7686954244834561] 
\draw   (890,150) -- (970,150) -- (970,198) -- (890,198) -- cycle ;
% Text Node
\draw (900,165) node [anchor=north west][inner sep=0.75pt]    {Loss};

%Straight Lines [id:da06856377424905746] 
\draw    (930,10) -- (930,140) ;
\draw [shift={(930,143)}, rotate = 270] [color={rgb, 255:red, 0; green, 0; blue, 0 }  ][line width=0.75]    (10.93,-3.29) .. controls (6.95,-1.4) and (3.31,-0.3) .. (0,0) .. controls (3.31,0.3) and (6.95,1.4) .. (10.93,3.29)   ;

%Straight Lines [id:da6090474583828176] 
\draw    (930,360) -- (930,203) ;
\draw [shift={(930,206)}, rotate = 90] [color={rgb, 255:red, 0; green, 0; blue, 0 }  ][line width=0.75]    (10.93,-3.29) .. controls (6.95,-1.4) and (3.31,-0.3) .. (0,0) .. controls (3.31,0.3) and (6.95,1.4) .. (10.93,3.29)   ;

\end{tikzpicture}
\end{center}
	\caption{Workflow  of the proposed PINN framework.}
	\label{fig:workflow}
\end{figure}

By following the penalization of the inequality constraints in the optimization theory \cite[Chapter~17]{JNocedal_SJWright_2006a}, we enforce the energy constraint  \eqref{energy_law}  to the loss function \eqref{lossEq} weakly by adding the following penalty term
\begin{equation}
	\hbox{loss}_e =  \frac{1}{n_e} \sum \limits_{k=1}^{n_e} \left( \max \left(0, \frac{\partial \mathcal{E}}{\partial t}\big(\widehat{u}(\bx_k^e, t_k^e)\big) \right)  \right)^2.
\end{equation}
Here, the energy functional approximation, denoted by $\mathcal{E}(\widehat{u})$, is obtained by evaluating  the approximated solution $\widehat{u}$ at the collocation points $(\bx_k^e, t_k^e)$, whereas the  derivative of the energy term with respect to time $t$ is obtained by using automatic differentiation as the previous computations. Finally, our total loss function becomes
\begin{equation}\label{loss_all}
\hbox{Loss}(\theta) = \lambda_r \, \hbox{loss}_r(\theta) + \lambda_i \, \hbox{loss}_i(\theta) +  \lambda_b \, \hbox{loss}_b(\theta) + \lambda_{e} \, \hbox{loss}_e(\theta),
\end{equation}
where $\lambda_r$, $\lambda_i$, $\lambda_b$, and $\lambda_e$ are the scale parameters of the domain residual, initial condition, boundary condition, and energy losses, respectively. Figure~\ref{fig:workflow} illustrates the framework of the proposed PINN approach to solve the Allen-Cahn equation.

%%%%%%%%%%%%%%%%%%%%%%%%%%%%%%%%%%%%%%%%%%%%%%%%%%%%%%%%%%%%%%%%%%%%%%%%%
\section{Adaptive approaches} \label{sec:algorithm} 

Now we provide some strategies based on adaptive approaches in the classical numerical techniques  to enhance the performance of the physics-informed neural network.

\subsection{Adaptive sampling of collocation points}

Phase field models as Allen-Cahn equation display sharp layers, where finer meshes are needed to capture the evolution dynamics in the traditional methods. Adaptive mesh refinement in the spatial domain has been considered as  one of the effective approaches in the literature to obtain a better accuracy with the lower degrees of freedom  for such kind of problems. Such a strategy  has been also applied successfully during the training processes in the several studies; see, e.g., \cite{LLu_XMeng_ZMao_GEKarniadakis_2021,CLWight_JZhao_2021,JMHanna_JVAguadoa_SCCardonab_RAskri_DBorzacchielloa_2022}. In a similar manner, we adaptively increase the number of collocation points by utilizing the residuals on the domain and boundary as an estimator; see Algorithm~\ref{alg:space}.
\begin{algorithm}[htp!]
	\caption{Adaptive Sampling of Collocation Points.}
	\label{alg:space}
	\begin{algorithmic}[1]
		\STATE Randomly generate  a new set of collocation points across the spatial domain \[\mathcal{V}:=\{(\bx_k^v, t_k^v)\}_{k=1}^{n_v},\] that is different from the current set of collocation points, denoted by $\mathcal{T}$.
		\STATE Compute the PDE residuals $\eta$ for the points in   $\mathcal{V}$.  %Make a prediction  on the data set $\mathcal{V}$   by applying a forward pass of the current neural network.
        \STATE Find a subset $\mathcal{M} \subset \mathcal{V}$ with the smallest cardinality satisfying the bulk criteria \cite{WDorfler_1996}:
        \[
           \tau  \sum \limits_{(\bx_k^{v_r},t_k^{v_r}) \in \mathcal{V}}    \eta(\bx_k^{v_r},t_k^{v_r})  \leq 
                  \sum \limits_{(\bx_k^{v_r},t_k^{v_r}) \in \mathcal{M} \subset \mathcal{V}}  \eta (\bx_k^{v_r},t_k^{v_r}),   
        \]
        where  $\tau \in (0,1)$ is a bulk parameter.
		\STATE Add the points in the set $\mathcal{M}$ into the set $\mathcal{T}$, that is, $\mathcal{T} \leftarrow \mathcal{T} \cup \mathcal{M}$.
		\STATE Continue training with the updated collocation set $\mathcal{T}$.
	\end{algorithmic}
\end{algorithm}
This algorithm not only allows us to select  samples from the sharp layers, but also to balance the computational cost by adaptively increasing the number of samples. On the other hand,  due to the computations in the second line of the algorithm, we run the algorithm if   the maximum of difference of loss values between last three successive iterations is greater than the given threshold $tol_{s}$, that is,
\begin{equation}
          \max(L_k, L_{k-1}, L_{k-2}) \geq tol_s  \quad \hbox{with} \quad  L_k = \hbox{Loss}(\theta^{k}) - \hbox{Loss}(\theta^{k-1}),
\end{equation}
and number of executed epochs is greater than a level denoted by $n_{ex}=1000$. The size of new data, $n_v$, is also taken as \%20 of the current training points.

\subsection{Adaptivity in the temporal domain}

As done in the spatial domain, it is also possible to take advantages of adaptivity in the temporal domain;  see, e.g., \cite{RMattey_SGhosh_2022,CLWight_JZhao_2021}. In this setting, separate networks for each time interval with the given time step $\Delta t$ are constructed. For the first network, we use the given initial condition $H(\bx)$. For the next time intervals, the predicted solution $\widehat{u}$ in the previous time segment is used as the initial condition. Algorithm~\ref{alg:time} summarizes the adaptivity process in the temporal domain.

%This process is repeated until the  value of energy functional becomes stable, that is,
%\begin{equation}\label{eq:energy_cond}
%\mathcal{E}(\widehat{u}(t_j)) - \mathcal{E}(\widehat{u}(t_{j-1})) \leq tol_a \quad \hbox{or} \quad \mathcal{E}(\widehat{u}(t_j)) \leq tol_b,
%\end{equation}
%where $tol_a, tol_b$ are the given threshold tolerance values. 

\begin{algorithm}[htp!]
	\caption{Adaptive Time Strategy.}
	\label{alg:time}
	\begin{algorithmic}[1]
		\STATE Given the time step $\Delta t$ and  maximum iteration number $N_{max}$. %and tolerance values $tol_a, tol_b$.
        \STATE Set $t_0= 0$.
        %\STATE Compute the initial energy $\mathcal{E}(\widehat{u}(t_0))$.
		\FOR{$j=1 \ldots N_{max}$}
		\STATE Set $t_j = t_{j-1} + \Delta t$ and  $I_j = [t_{j-1}, t_{j}]$.
        \IF {$(j=1)$}
        \STATE Set $H(\bx)$ as an initial condition for the interval $I_{j}$.
        \ELSE   
        \STATE Set the solution $\widehat{u}(\bx ,t_{j-1})$ as an initial condition for the interval $I_{j}$.
        \ENDIF
		\STATE Create a neural network on $I_j$ and predict the solution $\widehat{u}(x,t_{j})$  by using the trained neural network.
        %\IF{($\mathcal{E}(\widehat{u}(t_j)) - \mathcal{E}(\widehat{u}(t_{j-1})) \leq tol_a \lor \mathcal{E}(\widehat{u}(t_j)) \leq tol_b$)}
        %\STATE \textbf{break}.
        %\ENDIF
		\ENDFOR
	\end{algorithmic}
\end{algorithm}

\section{Numerical experiments}\label{sec:numeric} 

In this section, we will examine the effectiveness of the proposed energy dissipation approach outlined in Section~\ref{sec:energy_pinn} and compare it with existing studies in the literature.

Unless otherwise states, the hyperparameters used in the simulations are described in Table~\ref{tab:para}. Hyperparameters are key settings that significantly affect a model's performance. Consequently, hyperparameter optimization, which involves identifying the optimal combination of these settings to either minimize or maximize an objective function such as training loss or model accuracy, is essential. This optimization can be achieved through automated techniques like OPTUNA \cite{TAkiba_SSano_TYanase_TOhta_MKoyama_2019}, grid search \cite{PRen_CRao_YLiu_JXWang_HSun_2022}, Bayesian optimization \cite{YWang_XHan_CYChang_DZha_UBragaNeto_XHu_2023}, or
by employing a manual tuning strategy \cite{PPantidis_HEldababy_CMTagle_MEMobasher_2023,NHasebrook_FMorsbach_NKannengieser_MZoller_JFranke_MLindauer_FHutter_ASunyaev_2023}.  
However, in this study, when deciding some of the  hyperparameters, we consider the physical properties of the underlying partial differential equation (PDE), including its nonlinearity, time dependence, and the one-dimensional analytical solution for the transition layer, which describes the shift between the stable states 
$u = \pm 1$, that is,  $u= \tanh \left( \frac{x-x_0}{\sqrt{2} \epsilon}\right)$ for any arbitrary point $x_0$ in the domain. We prefer to use ResNet network structure \cite{KHe_2016} due to  its advantages in combining low-level and high-level features in the learning process.  To take advantages of both ADAM (first-order method) and L-BFGS (second-order method), we here apply ADAM+L-BFGS strategy as optimizer  \cite{LLu_XMeng_ZMao_GEKarniadakis_2021}, which  means that ADAM is applied until a certain number of iterations and then switched to L-BFGS method. As suggested in  \cite{DPKingma_JBa_2017}, the learning rate of ADAM optimizer is $0.001$ with all other parameters. With the help of this hybrid approach, ADAM can handle the initial exploration and rapid convergence, while L-BFGS can refine the solution by leveraging its second-order information. In the test examples, we have two pure states for the solution and the solution approaches one of them over time. Therefore, the hyperbolic tangent $\tanh$  is taken as an activation function  by considering its success in the two-class classification problems and the behaviour of the analytical solution for the transition layer. The activation function of the output layer, that is,  $\sigma_L$,  is a affine transformation for polynomial potentials, while for logarithmic potentials, it is the hyperbolic tangent function $\tanh$. To sample training points, we use  Latin hypercube sampling  \cite{MDShields_JZhang_2016}  for all cases. Further, all the neural networks are trained on Intel Core i7-6700 CPU 3.40GHz and 16 GB RAM. The software packages used for all simulations are Tensorflow 1.14 and MATLAB R2021a, and the data type of all the variables are float32. We note that the computations have been conducted on Central Processing Units (CPUs) due to resource constraints. However, for machine learning (ML) tasks, utilizing Graphical Processing Units (GPUs) can be more advantageous since GPUs are specifically designed to break down complex problems into thousands of smaller tasks and process them simultaneously.

%For finding the (optimal) hyperparameters in the PINNs applications, several methodologies have been suggested, such as  and manual tuning\cite{PPantidis_HEldababy_CMTagle_MEMobasher_2023,PSharma_LEvans_MTindall_PNithiarasu_2023,}. 

\begin{table}[htp!]
	\caption{Hyperparameters used in the  simulations.}\label{tab:para}
	\centering
    \begin{tabular}{ll}
	Hyperparameters & Used approaches                                             \\ \hline
	Network structure             & residual neural network (ResNet) \cite{KHe_2016}           \\
	Initialization of weights     & Xavier initialization \cite{XGlorat_Y.Bengio_2010}         \\
	Sampling of points            & Latin hypercube sampling  \cite{MDShields_JZhang_2016}     \\
	Activation function           & hyperbolic tangent $\tanh$                                 \\
    Optimizer                    & ADAM + L-BFGS  \cite{LLu_XMeng_ZMao_GEKarniadakis_2021}      \\
    \end{tabular}
\end{table}

In the numerical simulations, we begin by  comparing  the performance of our method with the existing PINN approaches from the literature using a well-studied example. 
To test accuracy of the obtained solution on the testing points, we use a relative total error over  the entire domain as done in \cite{MRaissi_PPerdikaris_GEKarniadakis_2019,CLWight_JZhao_2021,RMattey_SGhosh_2022}
\[
error = \frac{\left( \frac{1}{n} \sum_{i=1}^n |u_r(\bx_i,t_i) - \widehat{u}(\bx_i,t_i)|^2 \right)^{1/2}}{\left( \frac{1}{n} \sum_{i=1}^n |u_r(\bx_i,t_i)|^2 \right)^{1/2}},
\]
where $u_r$ is the analytic solution (if exists); otherwise,  it is a reference solution, and $\{(\bx_i, t_i)\}_{i=1}^n$ are test data points, which correspond to the grid points here. In the rest of simulations, the proposed approach will be tested by using various forms of the Allen-Cahn equation to demonstrate that the predicted solutions align with existing literature and adhere to the maximum bound principle, particularly ensuring the energy dissipation.

%%%%%%%%%%%%%%%%%%%%%%%%%%%%%%%%%%%%%%%%%%%%%%%%%%%%%%%%%%%%%%%%%%%%%%%%%%%%%%%%%%%%%%%%%%%%%%%%%%%%%%

\subsection{1D Allen--Cahn with polynomial potential}\label{ex:1D_double}
As a first benchmark problem, we consider one-dimensional Allen-Cahn equation with the periodic boundary conditions and the quartic double-well potential energy \eqref{Polynomial_Energy}, formulated as follows
\begin{equation} \label{eqn:1dexample}
    \begin{split}
	u_t - 0.0001u_{\bx\bx} + 5(u^3 - u) &= 0, \qquad \qquad \bx \, \in \, [-1, 1], \quad t \in [0,T], \\
	u(\bx, 0) &= \bx^2 \cos(\pi \bx), \\
	u(-1,t) &= u(1,t), \\
	u_\bx(-1,t) &= u_\bx(1,t).
    \end{split}
\end{equation}

In this example, we aim to compare the performance of our methodology with various PINN variants used to solve the benchmark problem in \eqref{eqn:1dexample}.
As done in the studies \cite{MRaissi_PPerdikaris_GEKarniadakis_2019,CLWight_JZhao_2021,RMattey_SGhosh_2022}, a reference solution  is obtained by using spectral methods with the  Chebfun package \cite{TADriscoll_NHale_LNTrefethen_2014}. To ensure a fair comparison with existing studies in the literature, hyperparameters are kept as consistent as possible; see, Table~\ref{tab:prm_ex1D}  for the used hyperparameters in this example.

\begin{table}[htp!]
	\caption{Description of training data  for the  Example \ref{ex:1D_double}.}\label{tab:prm_ex1D}
	\centering
    \resizebox{\textwidth}{!}
	{%
	\begin{tabular}{lll}
		Parameters                                 &    Values             & Descriptions                                             \\ \hline
		NN depth                                   &    4                  & \# hidden layer  \\
		NN width                                   &    100                & \# neurons in each hidden layer   \\
		$(n_r, n_b, n_i)$                          &  $(500,42,128)$       & size of initial samples per time segment       \\
		$(\lambda_r, \lambda_b, \lambda_i, \lambda_e)$ &  $(1,1,100,542)$  & scale parameters in \eqref{loss_all}       \\
		$(\tau, tol_s)$                            &  $(0.1,5\times 10^{-2})$         & parameters used in Algorithm~\ref{alg:space} \\
		$(\Delta_t, N_{max})$                      &  $(0.1,20)$           & parameters used in Algorithm~\ref{alg:time} \\
		$(N_{Adam},N_{LBFGS})$               	   &  $(5000,5000)$        & maximum ADAM and L-BFGS iterations per time segment \\
	\end{tabular}
    }
\end{table}

\begin{table}[htp!]
	\caption{Performance comparison for the  Example~\ref{ex:1D_double}.}
	\label{tab:results_ex1D}
	\centering
		\begin{tabular}{l|c}
			\textbf{method} & \textbf{relative $\ell_2$-norm}   \\
			\hline
			\texttt{std-PINN} \cite{MRaissi_PPerdikaris_GEKarniadakis_2019}   & 0.9919  \\
			\texttt{XPINN}    \cite{ADJagtap_GEKarniadakis_2020}              & 0.9612  \\
			\texttt{bc-PINN}  \cite{RMattey_SGhosh_2022}                      & 0.0701  \\
			\texttt{bc-PINN + logresidual} \cite{RMattey_SGhosh_2022}         & 0.0300  \\
			\texttt{ACP} \cite{CLWight_JZhao_2021}                            & 0.0233  \\
			\texttt{bc-PINN + ICGL + TL} \cite{RMattey_SGhosh_2022}           & 0.0168  \\
			\texttt{AT + ACP}                                                 & 0.0078 \\
			\texttt{AT + ACP + Energy Penalty}                                & 0.0053 
	\end{tabular}
\end{table}

\begin{figure}[htp!]
	\centering
	\includegraphics[width= 1\textwidth]{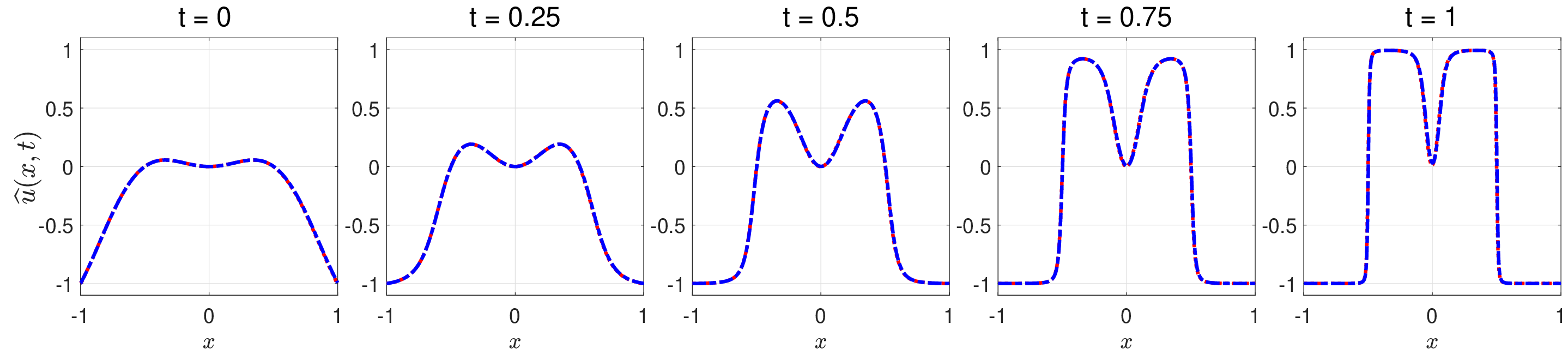}
    \includegraphics[width= 1\textwidth]{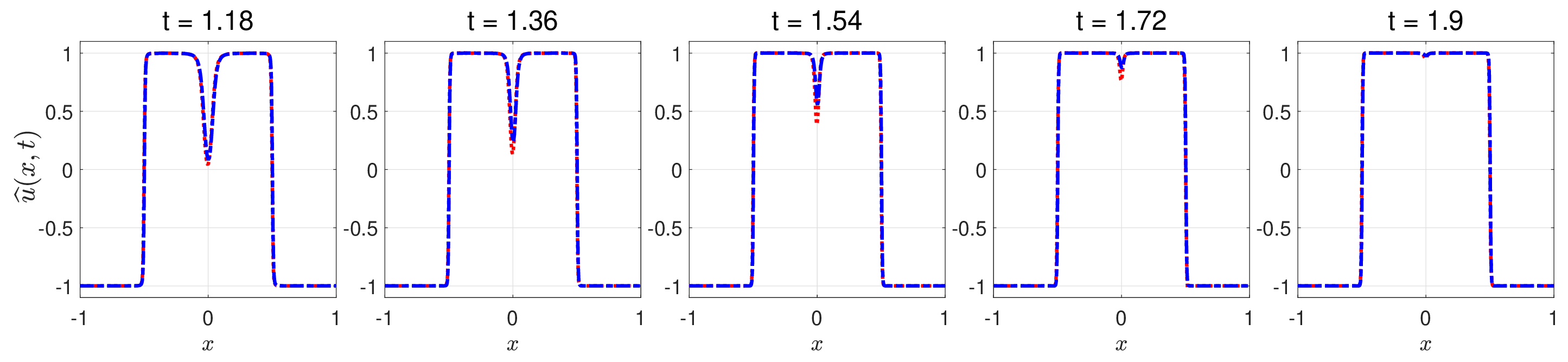}
	\vspace{-4mm} \caption{Example~\ref{ex:1D_double}: Predicted $\widehat{u}(x,t)$ (red) and reference $u_r(x,t)$ (blue) solutions at various time slots.}
	\label{fig:1D_ex1_soln}
\end{figure}  

\begin{figure}[htp!]
	\centering
	\includegraphics[width= 1\textwidth]{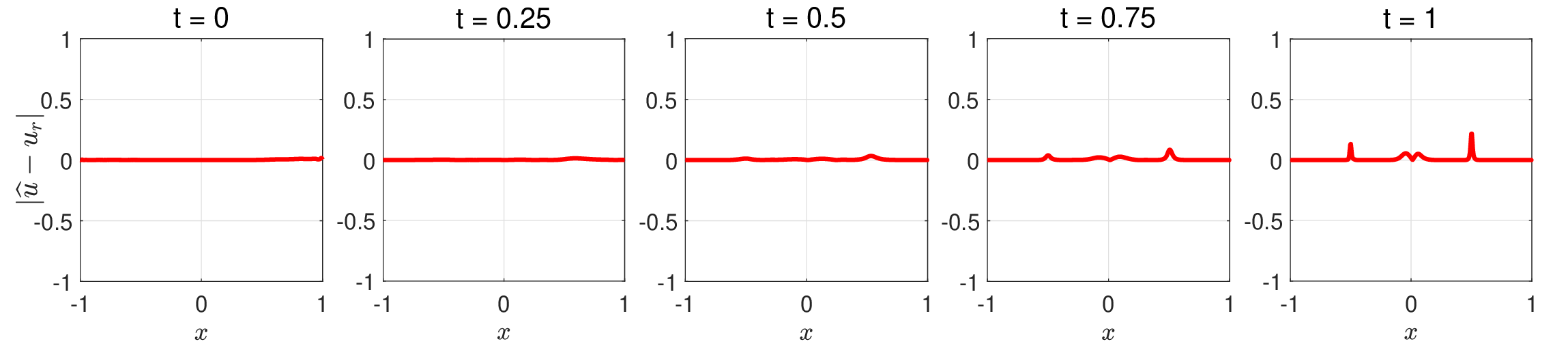}
    \includegraphics[width= 1\textwidth]{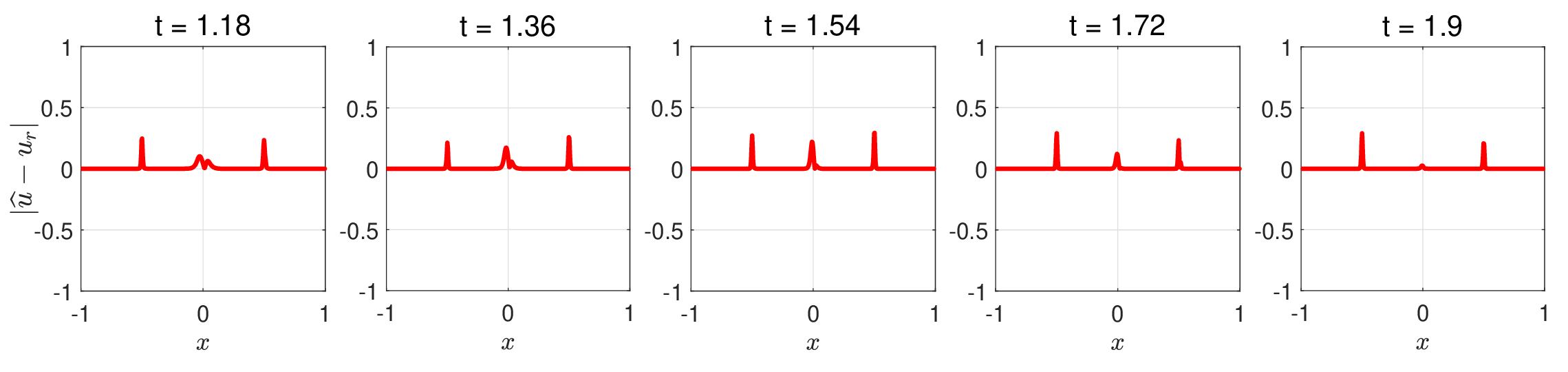}
	\vspace{-4mm} \caption{Example~\ref{ex:1D_double}: Behaviour of absolute error $|\widehat{u} - u_r |$ at various time slots.}
	\label{fig:1D_ex1_error}
\end{figure} 

In Table~\ref{tab:results_ex1D}, we compare the different approaches in the literature with our approach by taking the final time $T=1$.  By comparing with  \texttt{std-PINN}, \texttt{XPINN}, \texttt{bc-PINN}, \texttt{bc-PINN + logresidual}, and \texttt{bc-PINN + ICGL + TL} produced in \cite{RMattey_SGhosh_2022} with 4 hidden layer, 200 neurons in each layer, 512 points for initial condition, 201 points for the boundary condition, 20,000 collocation points for the residual per time segment, and with \texttt{ACP} produced in \cite{CLWight_JZhao_2021}  with 4 hidden layer, 128 neurons in each layer, 512 points for initial condition, 200 points for the boundary condition, 10,000 collocation points for the residual, our proposed approach  makes a significant improvement accuracy in spite of lower sampling points (less than 14,000 sampling points in total).

\begin{figure}[htp!]
	\centering
	\includegraphics[width= 0.45\textwidth]{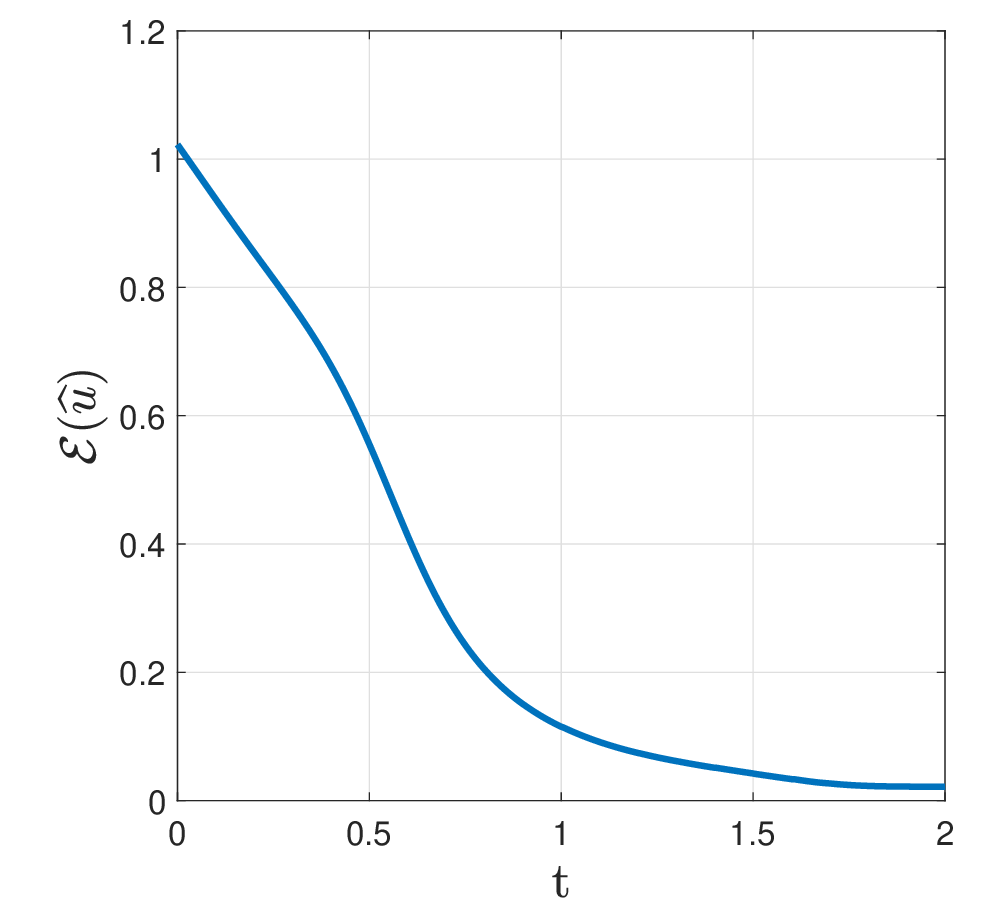}
	\caption{Example~\ref{ex:1D_double}: Time evolution of the predicted energy functional $\mathcal{E}(\widehat{u})$.}
	\label{fig:1D_ex1_energy}
\end{figure}

Figure~\ref{fig:1D_ex1_soln} presents a comparison between the predicted solution, $\widehat{u}$,  and the reference solution, $u_r$, obtained using Chebfun.
Although, the error increases with time slowly in Figure~\ref{fig:1D_ex1_error}, it is observed that their behaviours are almost similar.  In addition, our predicted solution stays bounded on [-1,1], which implies that the maximum bound principle is satisfied. Last, the behaviour of the predicted energy functional $\mathcal{E}(\widehat{u})$  is given in  Figure~\ref{fig:1D_ex1_energy}. As expected, it satisfies the decay of energy condition over time.

%%%%%%%%%%%%%%%%%%%%%%%%%%%%%%%%%%%%%%%%%%%%%%%%%%%%%%%%%%%%%%%%%%%%%%%%%%%%%%%%%%%%%%%%%%%%%%%%%%%%%%
\subsection{2D Allen--Cahn with polynomial potential}\label{ex:2D_double}
Next, we consider a two dimensional Allen-Cahn equation  with the quartic double-well potential energy \eqref{Polynomial_Energy}, in the form of 
\begin{eqnarray*}
	u_t &=& \mu(u) \big( \epsilon^2 \Delta u - u^3 + u \big), \hspace{41mm} \bx  \, \in \, \Omega, \quad t \in [0,T], \\
    \frac{\partial u}{\partial \mathbf{n}} &=& 0,  \hspace{77mm} \bx  \, \in \, \partial \Omega,  \\
	u(\bx, 0) &=& \tanh \left( \frac{0.35  - \sqrt{(x_1 - 0.5)^2 + (x_2 -0.5)^2}}{2 \epsilon} \right), \hspace{2mm}  \bx \in \Omega,
\end{eqnarray*}
where $\Omega= [0,1]^2$, $\epsilon=0.025$, and $\mu(u)=10$.  Hyperparameters used in this example are displayed in Table~\ref{tab:prm_ex2D}.

\begin{table}[htp!]
	\caption{Description of training data  for the  Example \ref{ex:2D_double}.}\label{tab:prm_ex2D}
	\centering
    \resizebox{\textwidth}{!}
	{%
	\begin{tabular}{lll}
		Parameters                                 &    Values                      & Descriptions                                             \\ \hline
		NN depth                                   &    6                           & \# hidden layer  \\
		NN width                                   &    128                         & \# neurons in each hidden layer   \\
		$(n_r, n_b, n_i)$                          &  $(5000,200,2601)$             & size of initial samples per time segment       \\
		$(\lambda_r, \lambda_b, \lambda_i, \lambda_e)$ &  $(1,1,1000,5200)$         & scale parameters in \eqref{loss_all}       \\
		$(\tau, tol_s)$                            &  $(0.05, 10^{-2})$             & parameters used in Algorithm~\ref{alg:space} \\
		$(\Delta_t, N_{max})$                      &  $(0.25,40)$                   & parameters used in Algorithm~\ref{alg:time} \\
		$(N_{Adam},N_{LBFGS})$               	   &  $(5000,2000)$                 & maximum ADAM and L-BFGS iterations per time segment  \\
	\end{tabular}
    }
\end{table}

\begin{figure}[htp!]
	\centering
	\includegraphics[width=0.9\textwidth]{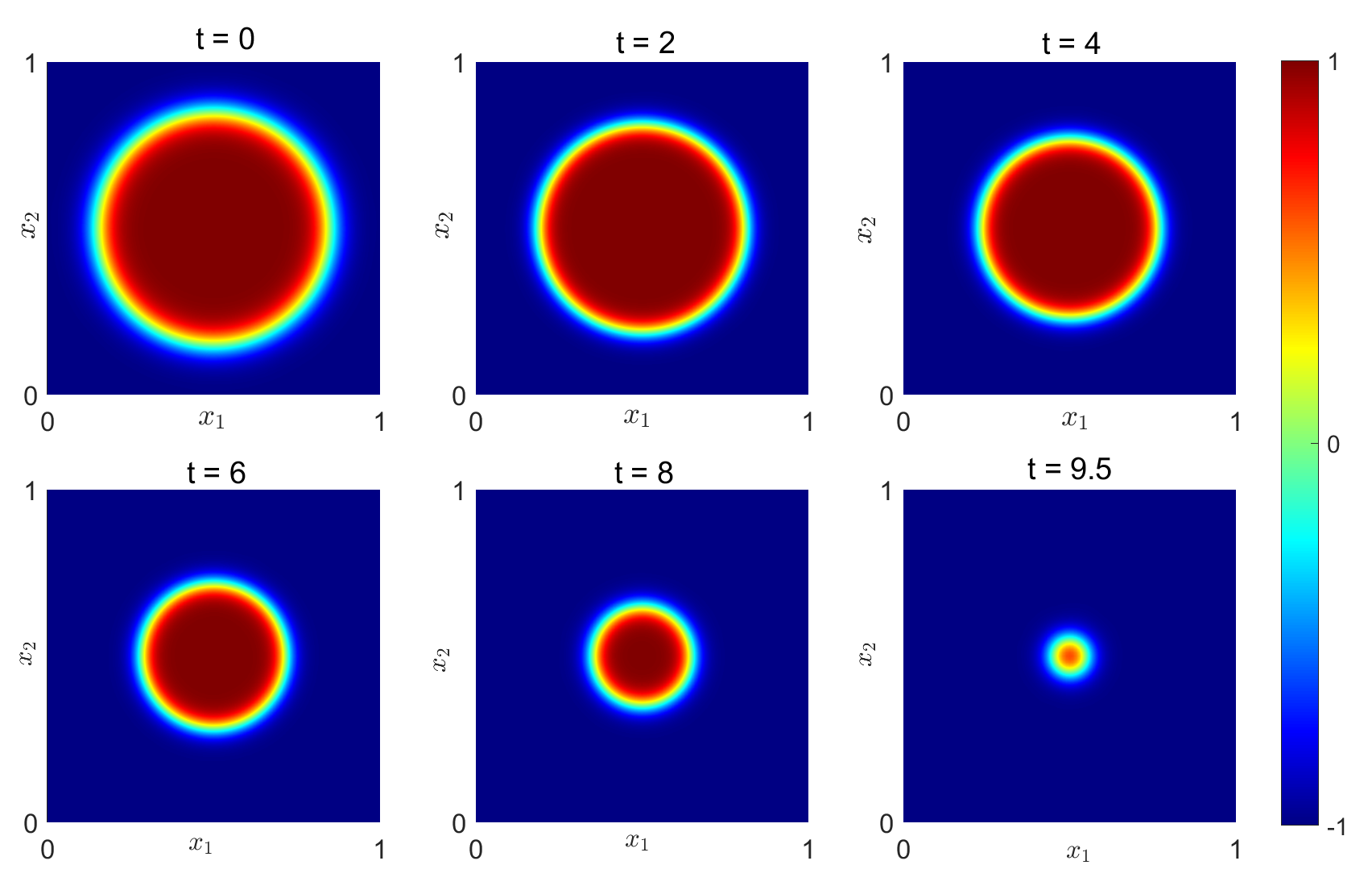}
	\caption{Example~\ref{ex:2D_double}: Predicted solutions $\widehat{u}(x,t)$  obtained by the proposed methodology at various time slots.}
	\label{fig:2D_ex2_soln}
\end{figure} 

\begin{figure}[htp!]
	\centering
	\includegraphics[width= 0.5\textwidth]{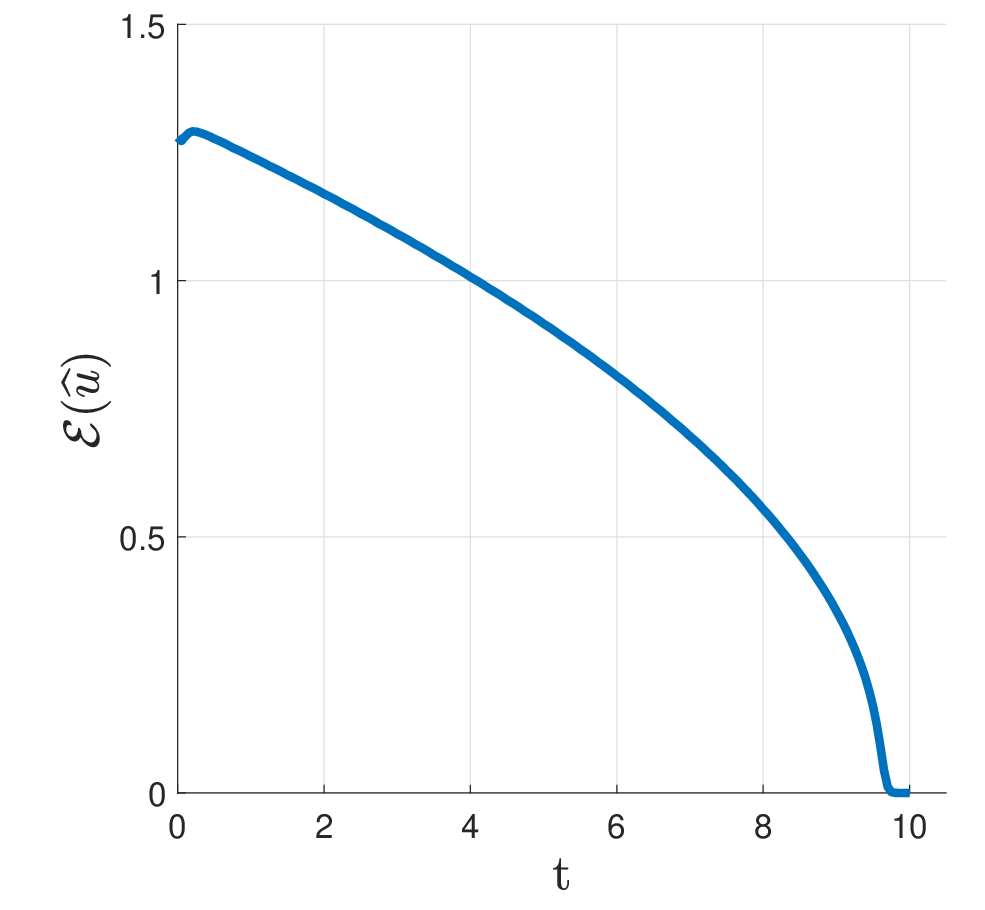}
	\caption{Example~\ref{ex:2D_double}: Time evolution of the predicted energy functional $\mathcal{E}(\widehat{u})$.}
	\label{fig:2D_ex2_energy}
\end{figure}

\begin{figure}[htp!]
	\centering
	\includegraphics[width=0.8\textwidth]{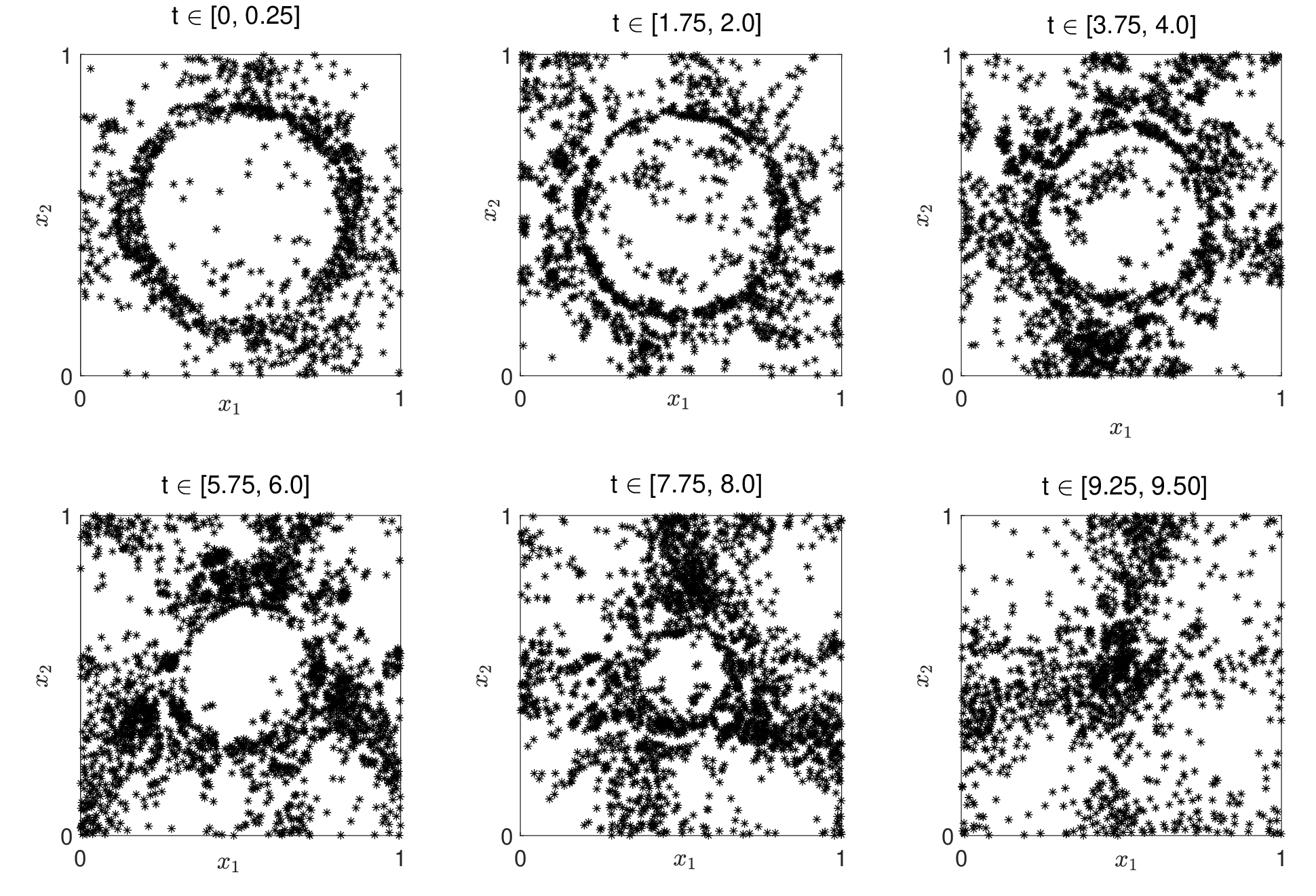}
	\vspace{-2mm} \caption{Example~\ref{ex:2D_double}: Distribution of adaptive collocation points at various time intervals.}
	\label{fig:ACP_Points_2D}
\end{figure}

Unlike the previous one-dimensional example, the trained weights and biases from the earlier time segment are transferred to the current time segment and serve as initials for the remaining test examples; see, e.g., \cite{SJPan_QYang_2010}. It is seen from Figure~\ref{fig:2D_ex2_soln} that the two-phrases, $u=-1$ and $u=1$, start moving and finally reaches the phrase $u=-1$. As time moves forward, the circles shrink with time to the center of the circle as expected; see, e.g., \cite{ASingh_RKSinha_2023,YLi_HGLee_DJeong_JKim}. The decay of the energy functional $\mathcal{E}(\widehat{u})$ over time is displayed in Figure~\ref{fig:2D_ex2_energy}. In Figure~\ref{fig:ACP_Points_2D},  we also investigate the performance the adaptive collocation approach in Algorithm~\ref{alg:space}. We notice that the adaptive collocation points are densely concentrated around the circle, which seems consistent with the results in Figure~\ref{fig:2D_ex2_soln}.

%%%%%%%%%%%%%%%%%%%%%%%%%%%%%%%%%%%%%%%%%%%%%%%%%%%%%%%%%%%%%%%%%%%%%%%%%%%%%%%%%%%%%%%%%%%%%%%%%%%%%%

\begin{figure}[htp!]
	\centering
	\includegraphics[width=0.9\textwidth]{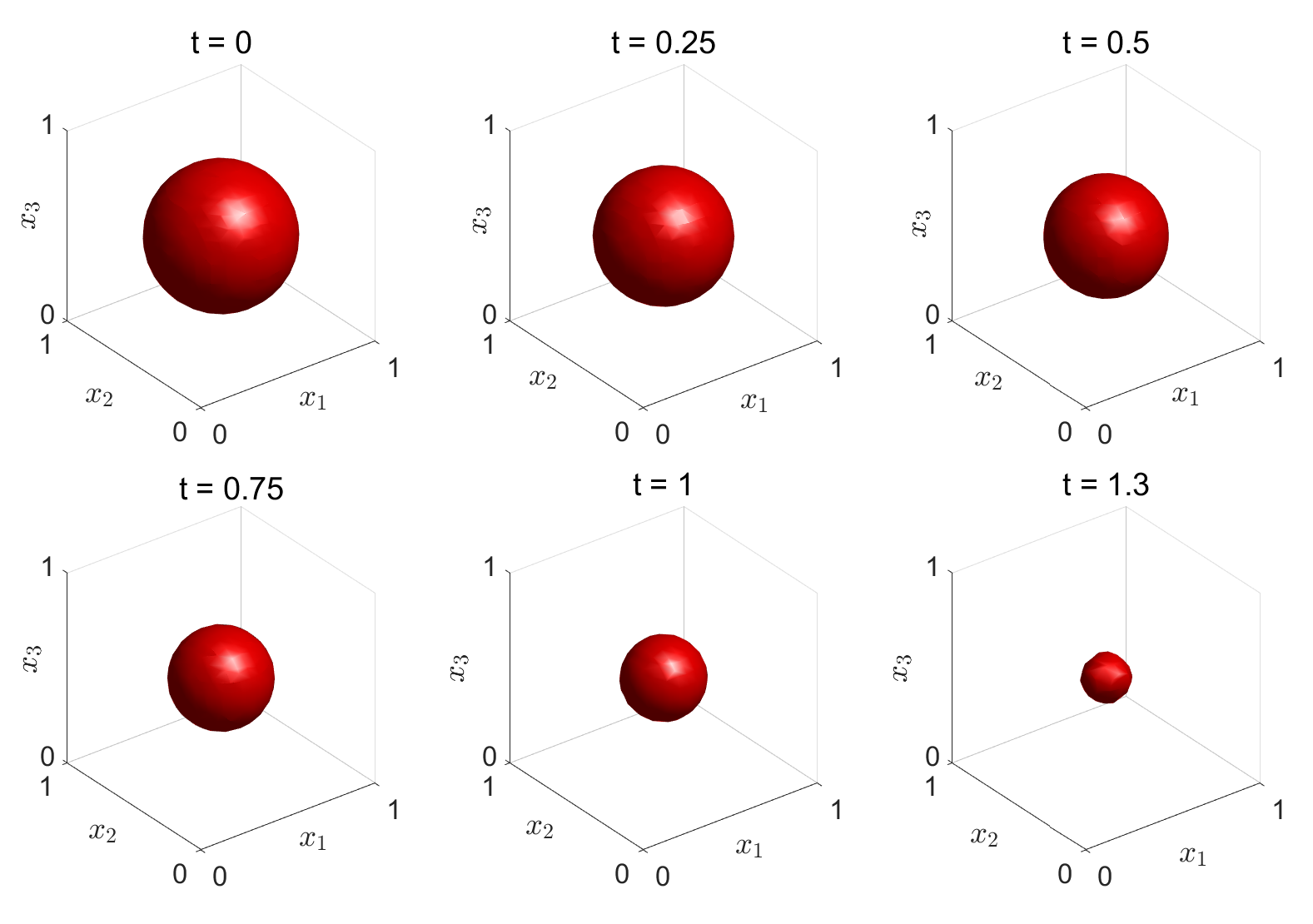}
	\caption{Example~\ref{ex:3D_double}: Predicted solutions $\widehat{u}(x,t)$ obtained by the proposed methodology at various time slots.}
	 \label{fig:3D_ex3_soln}
\end{figure}

\subsection{3D Allen--Cahn with polynomial potential}\label{ex:3D_double}
This example is an extension of the previous two dimensional Example \ref{ex:2D_double} into the domain $\Omega= [0,1]^3$ with the parameters  $\epsilon=0.05$,  $\mu(u)=10$, and the initial condition
\[
	u(\bx,0) = \tanh\Bigg(\frac{0.35 - \sqrt{(x_1-0.5)^2 + (x_2-0.5)^2 + (x_3-0.5)^2}}{2\epsilon}\Bigg).
\]

\begin{table}[htp!]
	\caption{Description of training data  for the  Example \ref{ex:3D_double}.}\label{tab:prm_ex3D}
	\centering
    \resizebox{\textwidth}{!}
	{%
	\begin{tabular}{lll}
		Parameters                                 &    Values                      & Descriptions                                             \\ \hline
		NN depth                                   &    6                           & \# hidden layer \\
		NN width                                   &    128                         & \# neurons in each hidden layer   \\
		$(n_r, n_b, n_i)$                          &  $(5000,300,9261)$             & size of initial samples per time segment       \\
		$(\lambda_r, \lambda_b, \lambda_i, \lambda_e)$ &  $(1,1,1000,5300)$         & scale parameters in \eqref{loss_all}      \\
		$(\tau, tol_s)$                            &  $(0.05, 10^{-2})$             & parameters used in Algorithm~\ref{alg:space} \\
		$(\Delta_t,N_{max})$                       &  $(0.1,13)$                    & parameters used in Algorithm~\ref{alg:time} \\
		$(N_{Adam},N_{LBFGS})$               	   &  $(5000,2000)$                 & maximum ADAM and L-BFGS iterations per time segment  \\
	\end{tabular}
    }
\end{table}

\begin{figure}[htp!]
	\centering
	\includegraphics[width=0.5\textwidth]{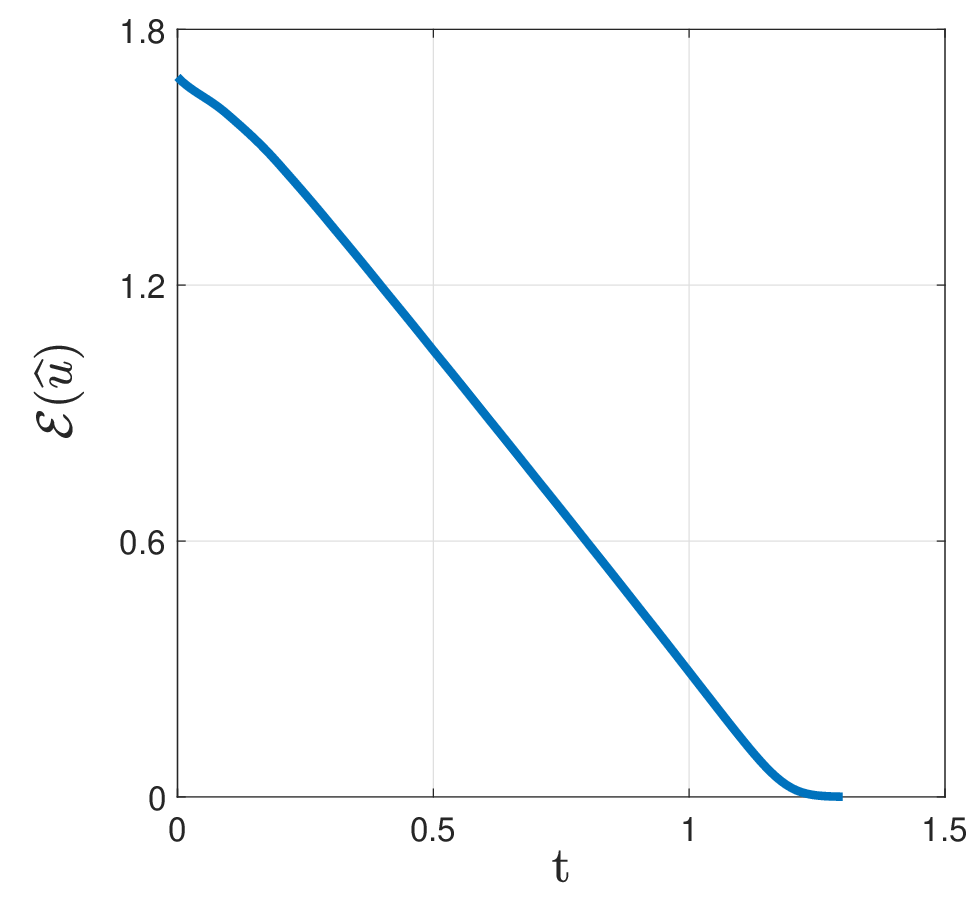}
	\caption{Example~\ref{ex:3D_double}: Time evolution of the predicted energy functional $\mathcal{E}(\widehat{u})$.}
	\label{fig:3D_ex3_energy}
\end{figure}

Table~\ref{tab:prm_ex3D} shows the used parameters in the simulation. Figure~\ref{fig:3D_ex3_soln} exhibits the zero isosurface of the solution at various time levels. Over time the size of the sphere decreases, which is also observed in \cite{ASingh_RKSinha_2023,YLi_HGLee_DJeong_JKim}. As expected, the predicted energy functional decreases $\mathcal{E}(\widehat{u})$ as time increases; see, Figure~\ref{fig:3D_ex3_energy}.

%%%%%%%%%%%%%%%%%%%%%%%%%%%%%%%%%%%%%%%%%%%%%%%%%%%%%%%%%%%%%%%%%%%%%%%%%%%%%%%%%%%%%%%%%%%%%%%%%%%%%%

\subsection{2D Allen--Cahn with logarithmic potential}\label{ex:2D_logorithmic}

Now, we consider a 2D Allen-Cahn equation \eqref{allencahn} with logarithmic potential \eqref{Flory_Huggins_Energy} on the domain $\Omega= [0, 2\pi]^2$, taken from \cite{DLi_CQuan_JXu_2022}. Rest of data is taken as 
\[
    \mu(u)=1, \quad \epsilon =0.1, \quad  \theta_c = 1, \quad  \theta=\frac{1}{4},
\]
and the initial condition consisting of seven circles
\begin{equation}\label{eq:2D_log_data}
    u(\bx,0) = -1 + \sum \limits_{i=1}^{7} h\Bigg(\sqrt{(x_1 - a_i)^2 + (x_2 - b_i)^2} -r_i \Bigg),
\end{equation}
where 
\[
h(s)= \left\{
       \begin{array}{ll}
         2 e^{-\epsilon^2/s^2}, & \hbox{if }  s<0, \\
         0, & \hbox{otherwise},
       \end{array}
     \right.
\]
and circles' centers and radii are given in Table~\ref{tab:2D_log_data}. Further, periodic boundary conditions are considered in this example; see~Table~\ref{tab:prm_ex2D_log} for the rest of used parameters in the simulation.

\begin{table}[htp!]
	\caption{Centers $(a_i,b_i)$ and radii $r_i$  in the initial condition \eqref{eq:2D_log_data} of the  Example~\ref{ex:2D_logorithmic}.}\label{tab:2D_log_data}
	\centering
	\begin{tabular}{c|lllllll}
		$i$   & $1$     & $2$       & $3$       & $4$       & $5$       & $6$      & $7$       \\ \hline
        $a_i$ & $\pi/2$ & $\pi/4$   & $\pi/2$   & $\pi$     & $3\pi/2$  & $\pi$    & $3\pi/2$  \\ 
        $b_i$ & $\pi/2$ & $3\pi/4$  & $5\pi/4$  & $\pi/4$   & $\pi/4$   & $\pi$    & $3\pi/2$  \\ 
        $r_i$ & $\pi/5$ & $2\pi/15$ & $2\pi/15$ & $\pi/10$  & $\pi/10$  & $\pi/4$  & $\pi/4$  \\ 
	\end{tabular}
\end{table}

\begin{table}[htp!]
	\caption{Description of training data  for the  Example~\ref{ex:2D_logorithmic}.}\label{tab:prm_ex2D_log}
	\centering
    \resizebox{\textwidth}{!}
	{%
	\begin{tabular}{lll}
		Parameters                                 &    Values                        & Descriptions                                             \\ \hline
		NN depth                                   &    6                             & \# hidden layer \\
		NN width                                   &    128                           & \# neurons in each hidden layer  \\
		$(n_r, n_b, n_i)$                          &  $(10000,200,40401)$             & size of initial samples per time segment       \\
		$(\lambda_r, \lambda_b, \lambda_i, \lambda_e)$ &  $(1,1,1000,10200)$          & scale parameters in \eqref{loss_all}       \\
		$(\tau, tol_s)$                            &  $(0.05, 10^{-1})$               & parameters used in Algorithm~\ref{alg:space} \\
		$(\Delta_t,N_{max})$                       &  $(0.125,40)$                    & parameters used in Algorithm~\ref{alg:time} \\
		$(N_{Adam},N_{LBFGS})$               	   &  $(5000,2000)$                   & maximum ADAM and L-BFGS iterations per time segment \\
	\end{tabular}
    }
\end{table}

\begin{figure}[htp!]
	\centering
	\includegraphics[width=1\textwidth]{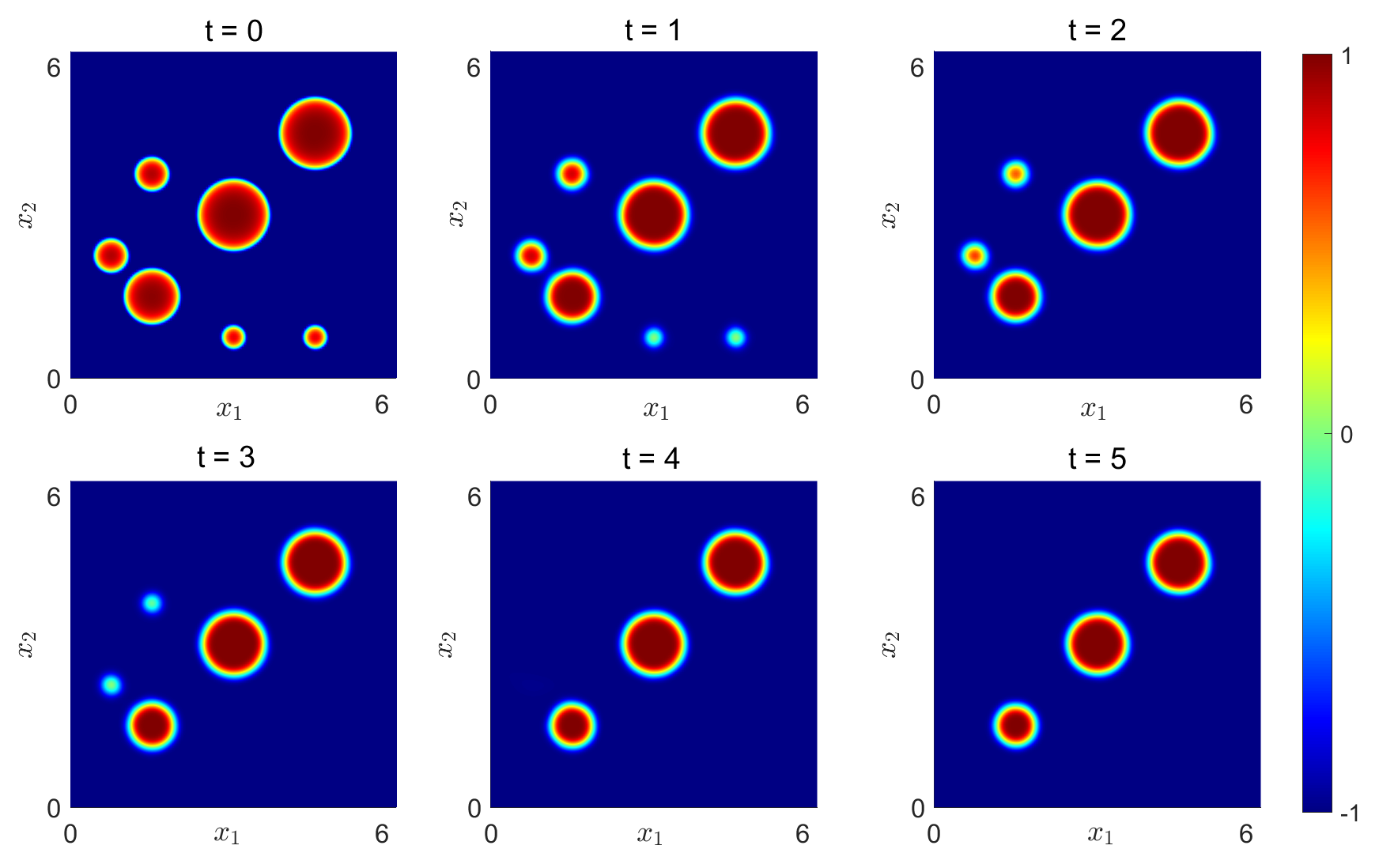}
	\caption{Example~\ref{ex:2D_logorithmic}: Predicted solutions $\widehat{u}(x,t)$ obtained by the proposed methodology at various time slots.}
	\label{fig:2D_log_ex4_soln}
\end{figure}  

\begin{figure}[h!]
	\centering
	\includegraphics[width=0.55\textwidth]{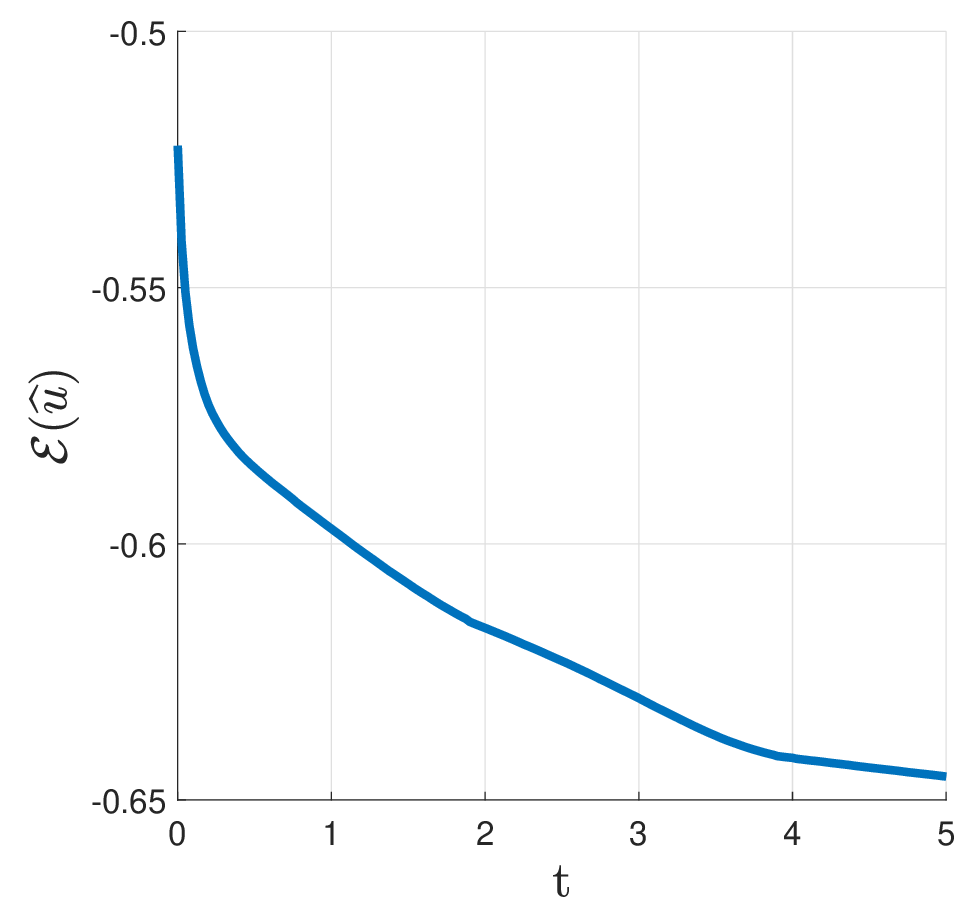}
	\caption{Example~\ref{ex:2D_logorithmic}: Time evolution of the predicted energy functional $\mathcal{E}(\widehat{u})$.}
	\label{fig:2D_log_ex4_energy}
\end{figure}

\begin{figure}[h!]
	\centering
	\includegraphics[width=0.8\textwidth]{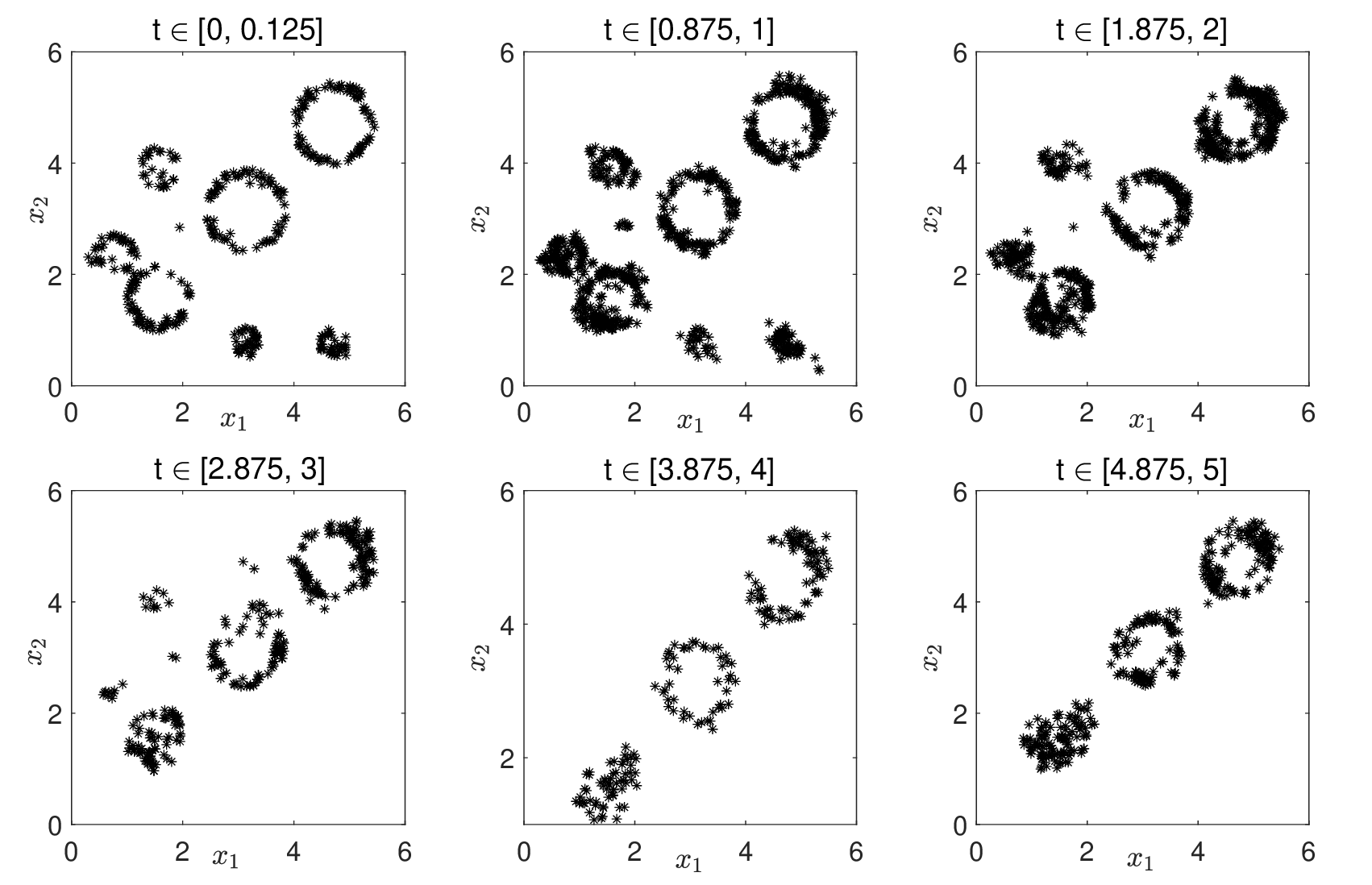}
	\caption{Example~\ref{ex:2D_logorithmic}: Distribution of adaptive collocation points at various time intervals.}
	\label{fig:ACP_Points_2D_log}
\end{figure}

From  the predicted solutions  in Figure~\ref{fig:2D_log_ex4_soln}, it is seen that the circles annihilate gradually in time as observed in \cite{DLi_CQuan_JXu_2022}. The behaviour of the approximated energy functional $\mathcal{E}(\widehat{u})$ is displayed in Figure~\ref{fig:2D_log_ex4_energy}. It seems that our approach results in a progressively diminishing estimated energy. Last, Figure~\ref{fig:ACP_Points_2D_log} shows that the distribution of adaptive collocation points become intense around the circles as expected from   the results in Figure~\ref{fig:2D_log_ex4_soln}.

%%%%%%%%%%%%%%%%%%%%%%%%%%%%%%%%%%%%%%%%%%%%%%%%%%%%%%%%%%%%%%%%%%%%%%%%%%%%%%%%%%%%%%%%%%%%%%%%%%%%%%

\subsection{2D Allen--Cahn with logarithmic potential and random initials}\label{ex:2D_logorithmic_random}

Next, we consider a benchmark example  \eqref{allencahn} with the logarithmic free energy  \eqref{Flory_Huggins_Energy}, periodic boundary conditions, and a random initial condition
\[
    u_0(x,0) = 0.05(2\times \hbox{\texttt{randn}} - 1)
\]
on the domain $\Omega= [0, 2\pi]^2$, taken from \cite{JShen_TTang_JYang_2016}. Here, \texttt{randn} generates normally distributed random vector. The rest of the problem data are
\[
  \mu(u)=2, \quad \epsilon = 0.04, \quad \theta=0.15, \quad \theta_c=0.30.
\]

Different from the previous examples, the learning process of the initial condition is not straightforward due to its inherent randomness. As a remedy, we use a continuous analogue of the initial random vector. For a complex periodic smooth random function on the domain $[-L/2, L/2]^2$ with $L>0$, its continuous analogue is characterized by a finite Fourier series with independent normally distributed random coefficients
\begin{equation}
	f(x,y) = \sum \limits_{k=-m}^{m} \sum \limits_{j=-m_k}^{m_k} c_{jk} \exp\Bigg(\frac{2 \pi \, i (jx+ky)}{L}\Bigg),
\end{equation}
where $m = \lfloor L/\gamma \rfloor$ with a wavelength parameter $\gamma>0$ and $m_k = \sqrt{m^2 - k^2}$. The coefficient parameters, $c_{jk}$, are  the element of a $(2m + 1) \times (2m_k + 1)$ matrix that includes independent sample from a normal distribution. We construct the corresponding smooth initial condition by using \texttt{randnfun2} built-in function in the  Chebfun package \cite{TADriscoll_NHale_LNTrefethen_2014}, regarded as a continuous analogue of the MATLAB command \texttt{randn}. We refer to \cite{SFilip_AJaveed_LNTrefethen_2019} for a more detailed discussion. In the simulations, the wavelength parameter $\gamma$ is taken as $\gamma=1$. It is not presented here but different values of $\gamma$ yields similar behaviours. 

\begin{table}[htp!]
	\caption{Description of training data  for the 
    Example~\ref{ex:2D_logorithmic_random}.}\label{tab:prm_ex2D_logorithmic_random} 
	\centering
    \resizebox{\textwidth}{!}
	{%
	\begin{tabular}{lll}
		Parameters                                 &    Values               & Descriptions                                             \\ \hline
		NN depth                                   &    6                    & \# hidden layer \\
		NN width                                   &    128                  & \# neurons in each hidden layer   \\
		$(n_r, n_b, n_i)$                          &  $(5000,200,40000)$     & size of initial samples per time segment       \\
		$(\lambda_r, \lambda_b, \lambda_i, \lambda_e)$ &  $(1,1,1000,5200)$      & scale parameters in \eqref{loss_all}       \\
		$(\tau, tol_s)$                            &  $(0.05,0.1)$            & parameters used in Algorithm~\ref{alg:space} \\
		$(\Delta_t,N_{max})$                       &  $(1,20)$                & parameters used in Algorithm~\ref{alg:time} \\
		$(N_{Adam},N_{LBFGS})$               	   &  $(5000,2000)$          & maximum ADAM and L-BFGS iterations per time segment \\
	\end{tabular}
    }
\end{table}

\begin{figure}[htp!]
 	\centering
 	\includegraphics[width=1\textwidth]{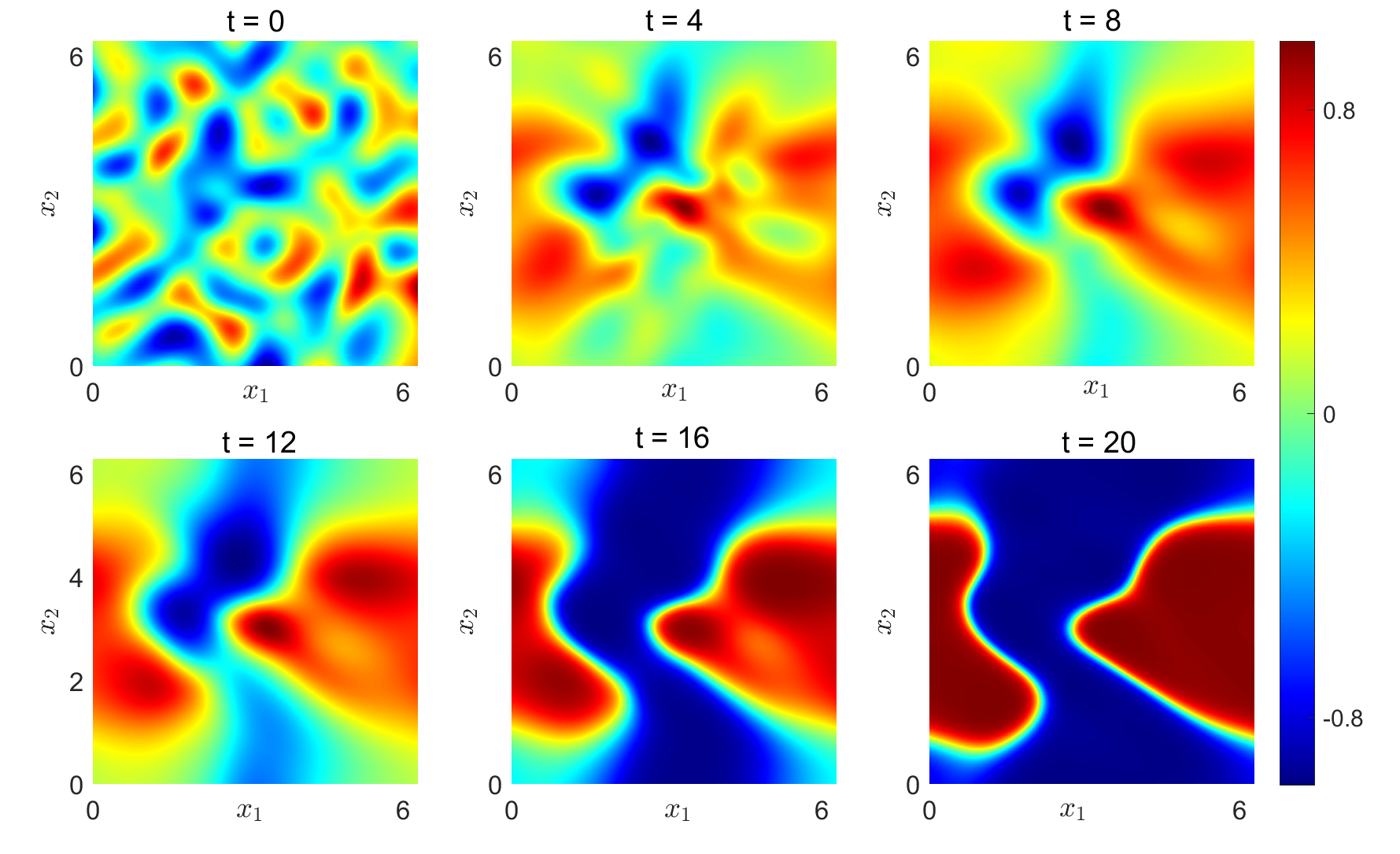}
 	\caption{Example~\ref{ex:2D_logorithmic_random}: Predicted solutions $\widehat{u}(x,t)$ obtained by the proposed methodology at various time slots.}
 	\label{fig:2D_log_rand_soln}
\end{figure}

\begin{figure}[htp!]
	\centering
	\includegraphics[width=0.5\textwidth]{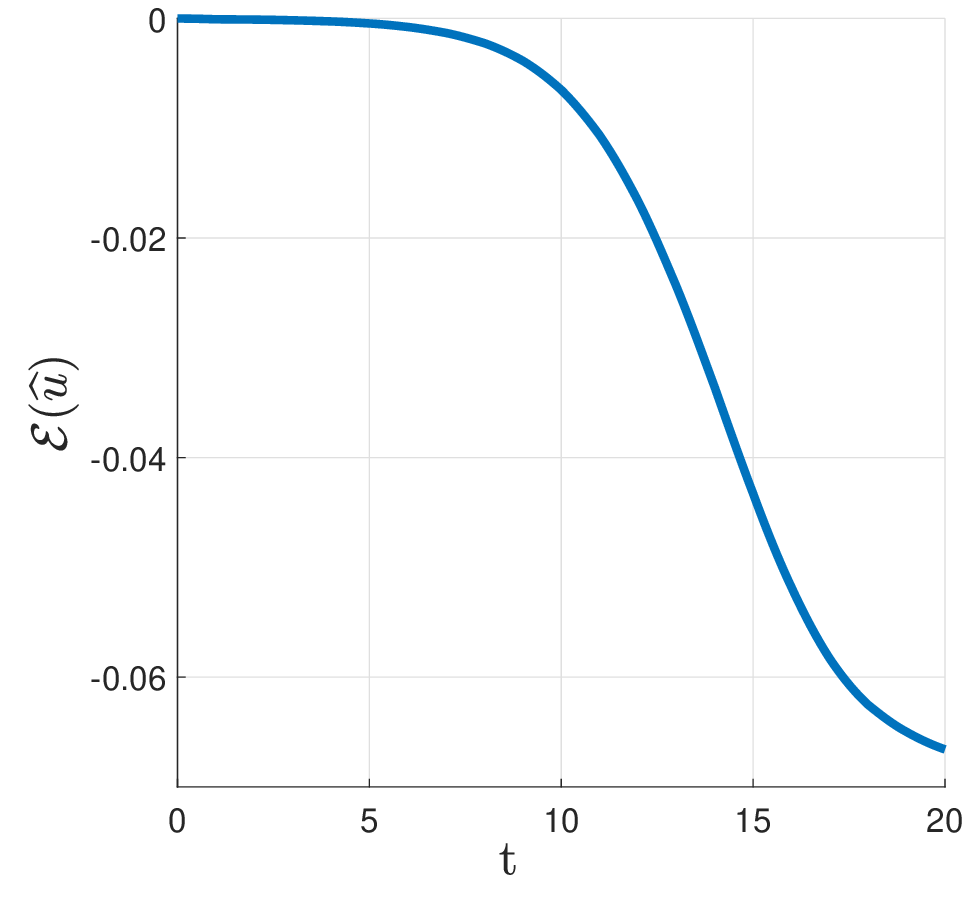}
	\caption{Example~\ref{ex:2D_logorithmic_random}: Time evolution of the predicted energy functional $\mathcal{E}(\widehat{u})$.}
	\label{fig:2D_log_rand_energy}
\end{figure}  

Table~\ref{tab:prm_ex2D_logorithmic_random} displays the parameters used in the training of the network. From Figure~\ref{fig:2D_log_rand_soln}, it is seen that the our approach yields a coarsening phenomenon during the phase evolution as expected; see also  \cite{JShen_TTang_JYang_2016}. In addition, the predicted solutions $\widehat{u}(x,t)$ is bounded in the interval $[-s,s]$  for all the time, where $1>s>0$ is the constant satisfying $f(s)=0$. It implies that the discrete maximum principle are indeed well satisfied. Further, the decay of energy function is displayed in Figure~\ref{fig:2D_log_rand_energy}.

%%%%%%%%%%%%%%%%%%%%%%%%%%%%%%%%%%%%%%%%%%%%%%%%%%%%%%%%%%%%%%%%%%%%%%%%%%%%%%%%%%%%%%%%%%%%%%%%%%%%%%

\begin{table}[htp!]
	\caption{Description of training data  for the  Example \ref{ex:2D_logorithmic_degenerate}.}\label{tab:prm_ex2D_logorithmic_degenerate} 
	\centering
    \resizebox{\textwidth}{!}
	{%
	\begin{tabular}{lll}
		Parameters                                 &    Values             & Descriptions                                             \\ \hline
		NN depth                                   &    6                  & \# hidden layer \\
		NN width                                   &    128                & \# neurons in each hidden layer   \\
		$(n_r, n_b, n_i)$                          &  $(7500,200,40.000)$  & size of initial samples per time segment       \\
		$(\lambda_r, \lambda_b, \lambda_i, \lambda_e)$ &  $(1,1,10^5,7700)$    & scale parameters in \eqref{loss_all}       \\
		$(\tau, tol_s)$                            &  $(0.05,0.1)$         & parameters used in Algorithm~\ref{alg:space} \\
		$(\Delta_t,N_{max})$                       &  $(1,10)$             & parameters used in Algorithm~\ref{alg:time} \\
		$(N_{Adam},N_{LBFGS})$               	   &  $(3000,20.000)$          & maximum ADAM and L-BFGS iterations per time segment \\
	\end{tabular}
    }
\end{table}

\subsection{2D Allen--Cahn with logarithmic potential and degenerate mobility}\label{ex:2D_logorithmic_degenerate}

Now we examine a two-dimensional Allen-Cahn equation characterized by the logarithmic potential \eqref{Flory_Huggins_Energy}, random initial condition, and degenerate mobility function within the domain $\Omega= [0, 2\pi]^2$. The functions and parameters of the underlying problem, taken from  \cite{BKarasozen_MUzunca_ASFilibelioglu_Hyucel_2018,JShen_TTang_JYang_2016}, are
\[
\epsilon = 0.04, \quad  \theta_c = 0.95, \quad  \theta = 0.5, \quad \mu = 2(1-u^2),
\]
with the random initial condition
\begin{equation*}
u(\bx,0) = 0.05(2 \times \hbox{\texttt{randn}} - 1).
\end{equation*}

%\begin{figure}[htp!]
% 	\centering
% 	\includegraphics[width=0.95\textwidth]{2D_AC_Log_Degenerate_IC_Comparison}
% 	\caption{Example~\ref{ex:2D_logorithmic_degenerate}: Comparison of the smooth initial condition with $\lambda=0.4$ and the predicted initial condition.}
% 	\label{fig:2D_log_degenerate_initial}
%\end{figure}  

\begin{figure}[t!]
 	\centering
 	\includegraphics[width=0.90\textwidth]{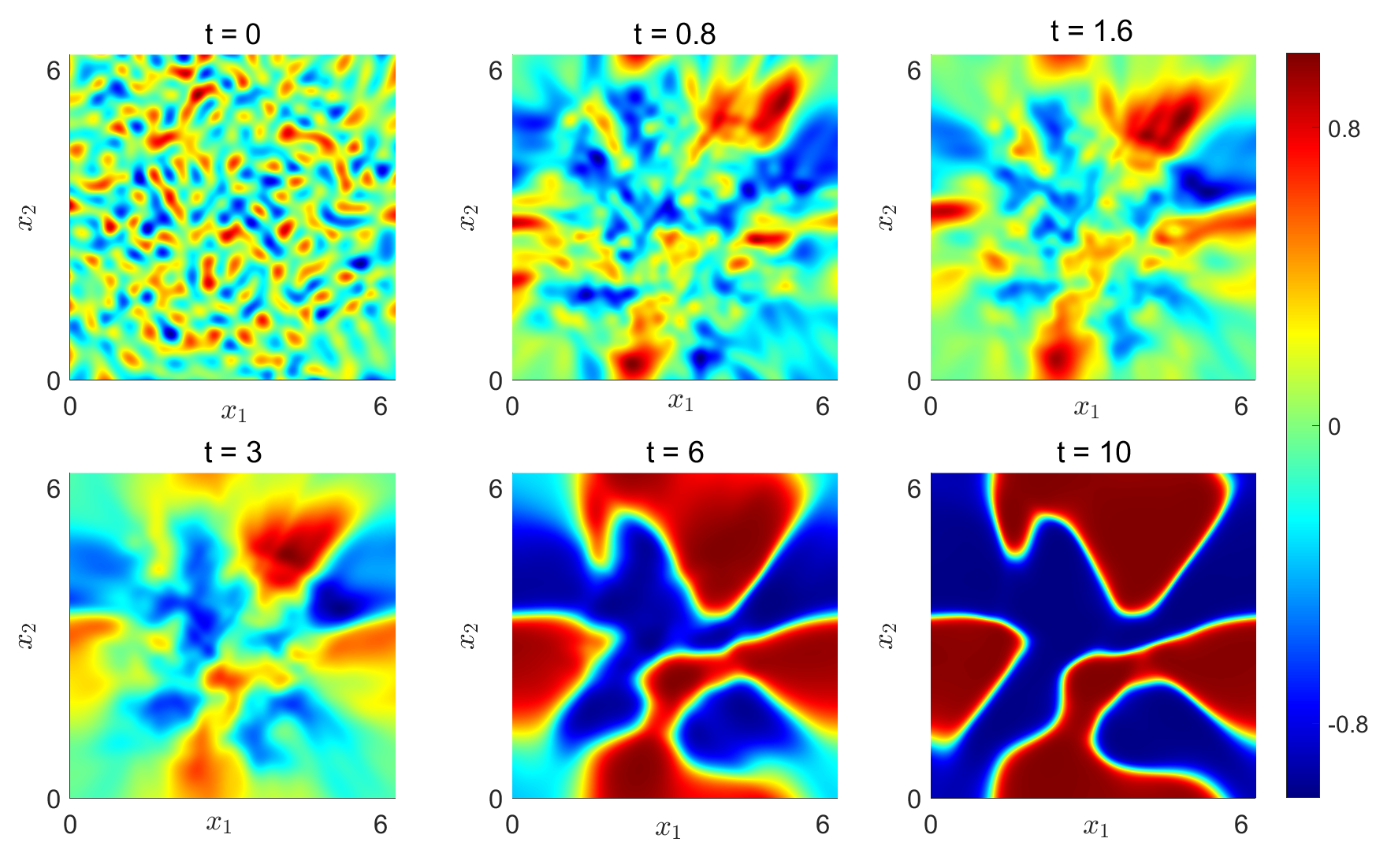}
 	\caption{Example~\ref{ex:2D_logorithmic_degenerate}: Predicted solutions $\widehat{u}(x,t)$ obtained by the proposed methodology at various time slots.}
 	\label{fig:2D_log_degenerate_soln}
\end{figure}  

\begin{figure}[htp!]
	\centering
	\includegraphics[width=0.38\textwidth]{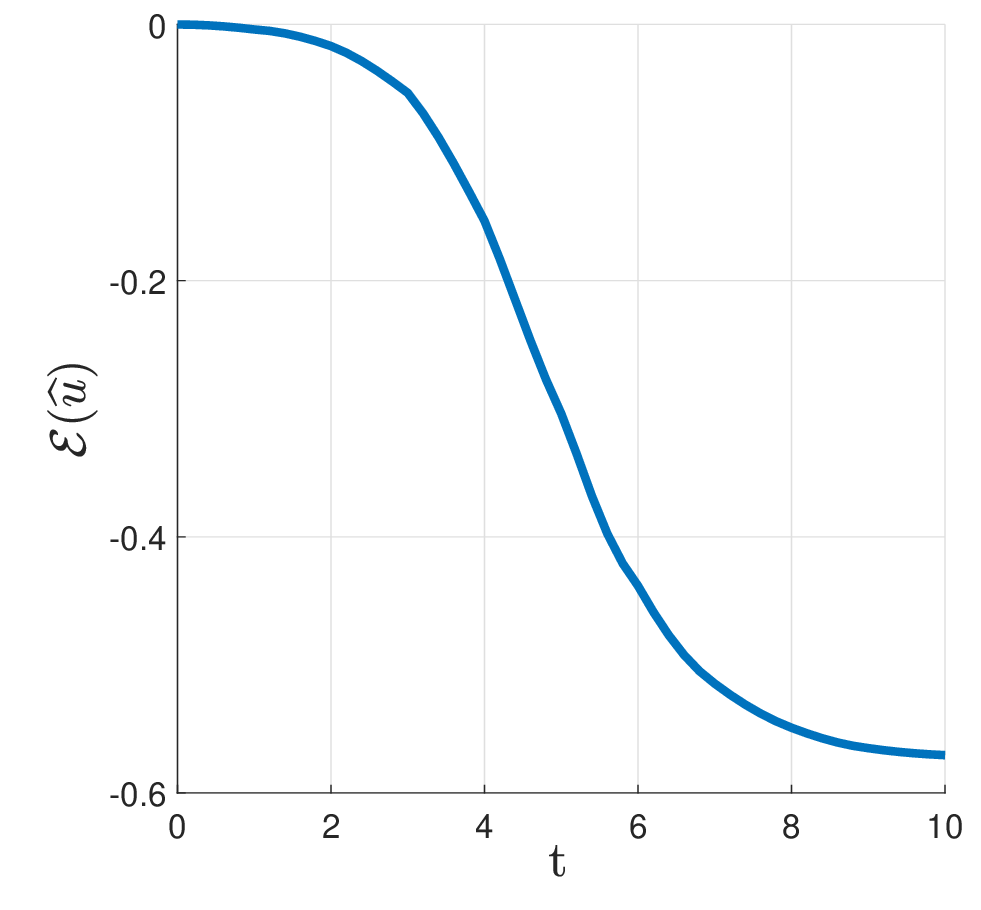}
	\caption{Example~\ref{ex:2D_logorithmic_degenerate}: Time evolution of the predicted energy functional $\mathcal{E}(\widehat{u})$.}
	\label{fig:2D_log_degenerate_energy}
\end{figure}  

As the previous example, we construct a continuous analogue of the random initial condition with the wavelength parameter $\gamma=0.4$. To address the challenge posed by the nonlinear mobility function, we adopt a strategy similar to explicit-implicit methods used in the numerical analysis. Instead of directly dealing with the nonlinear mobility function, we utilize the predicted  solutions from the previous time step in the Algorithm~\ref{alg:time}. Figure~\ref{fig:2D_log_degenerate_soln} illustrates the temporal behavior of the predicted solutions utilizing the parameters outlined in  Table~\ref{tab:prm_ex2D_logorithmic_degenerate}. The results show that a spinodal decomposition happens over time and the predicted solutions $\widehat{u}(x,t)$ preserve the boundedness. Figure~\ref{fig:2D_log_degenerate_energy} also displays that the approximate energy functional  $\mathcal{E}(\widehat{u})$ is dissipated.

%%%%%%%%%%%%%%%%%%%%%%%%%%%%%%%%%%%%%%%%%%%%%%%%%%%%%%%%%%%%%%%%%%%%%%%%%%%%%%%%%%%%%%%%%%%%%%%%%%%%%%

\begin{table}[H]
	\caption{Description of training data for the  Example~\ref{ex:2D_advection}.}\label{tab:prm_ex2D_advection} 
	\centering
    \resizebox{\textwidth}{!}
	{%
	\begin{tabular}{lll}
		Parameters                                 &    Values             & Descriptions            \\ \hline
		NN depth                                   &    6                  & \# hidden layer \\
		NN width                                   &    128                & \# neurons in each hidden layer   \\
		$(n_r, n_b, n_i)$                          &  $(10000,200,10201)$  & size of initial samples per time segment       \\
		$(\lambda_r, \lambda_b, \lambda_i, \lambda_e)$ &  $(1,1,10^{5},10200)$ & scale parameters in \eqref{loss_all}       \\
		$(\tau, tol_s)$                            &  $(0.05,0.1)$       & parameters used in Algorithm~\ref{alg:space} \\
		$(\Delta_t,N_{max})$                          &  $(0.005,12)$        & parameters used in Algorithm~\ref{alg:time} \\
		$(N_{Adam},N_{LBFGS})$               	   &  $(5000,2000)$                   & maximum ADAM and L-BFGS iterations per time segment \\
	\end{tabular}
}
\end{table}

\subsection{2D Allen--Cahn with advection and polynomial potential}\label{ex:2D_advection}

\begin{figure}[htp!]
	\centering
 	\includegraphics[width=1\textwidth]{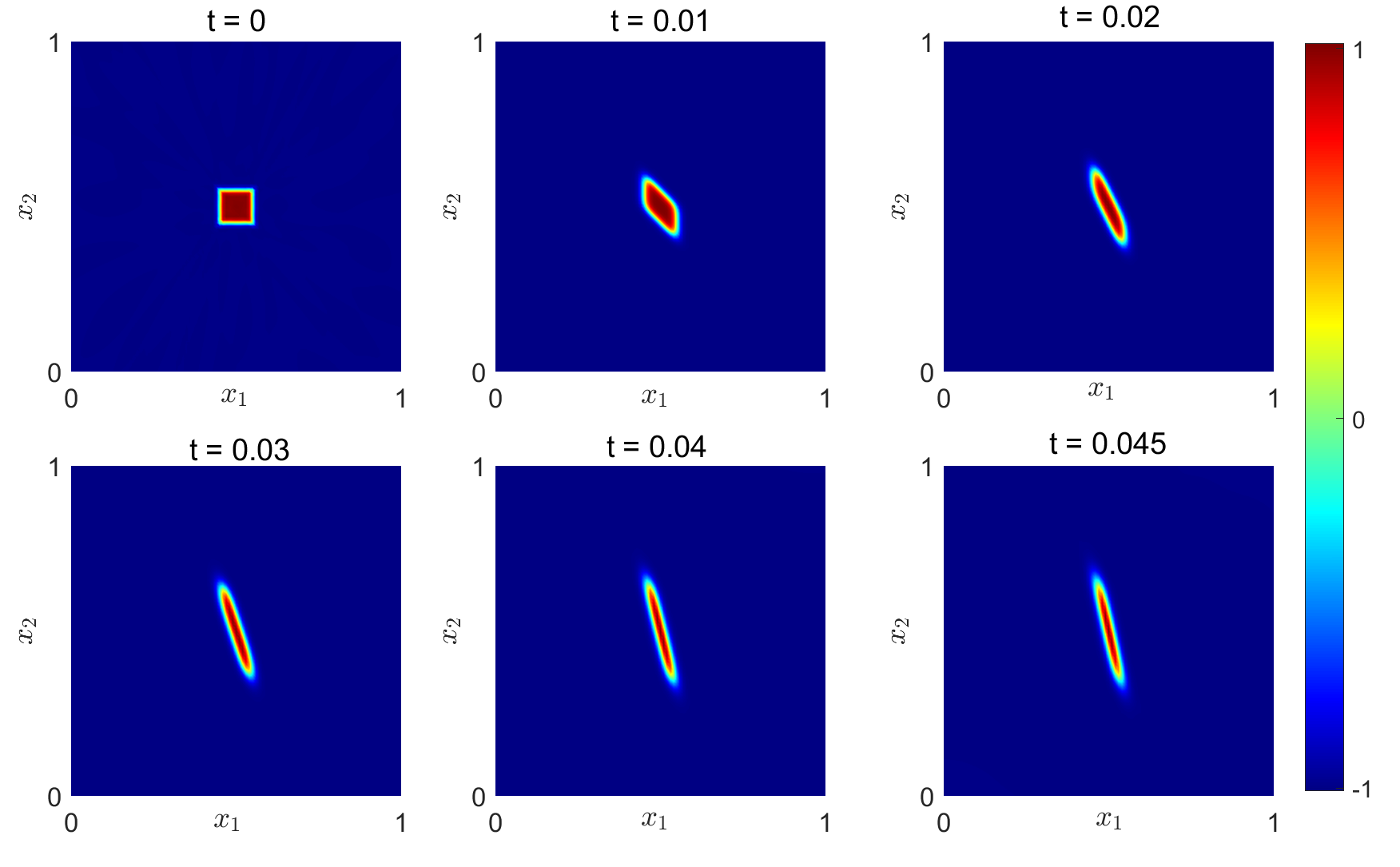}
 	\caption{Example~\ref{ex:2D_advection}: Predicted solutions $\widehat{u}(x,t)$ obtained by the proposed methodology at various time slots.}
 	\label{fig:2D_advection}
\end{figure}

\begin{figure}[htp!]
	\centering
	\includegraphics[width=0.5\textwidth]{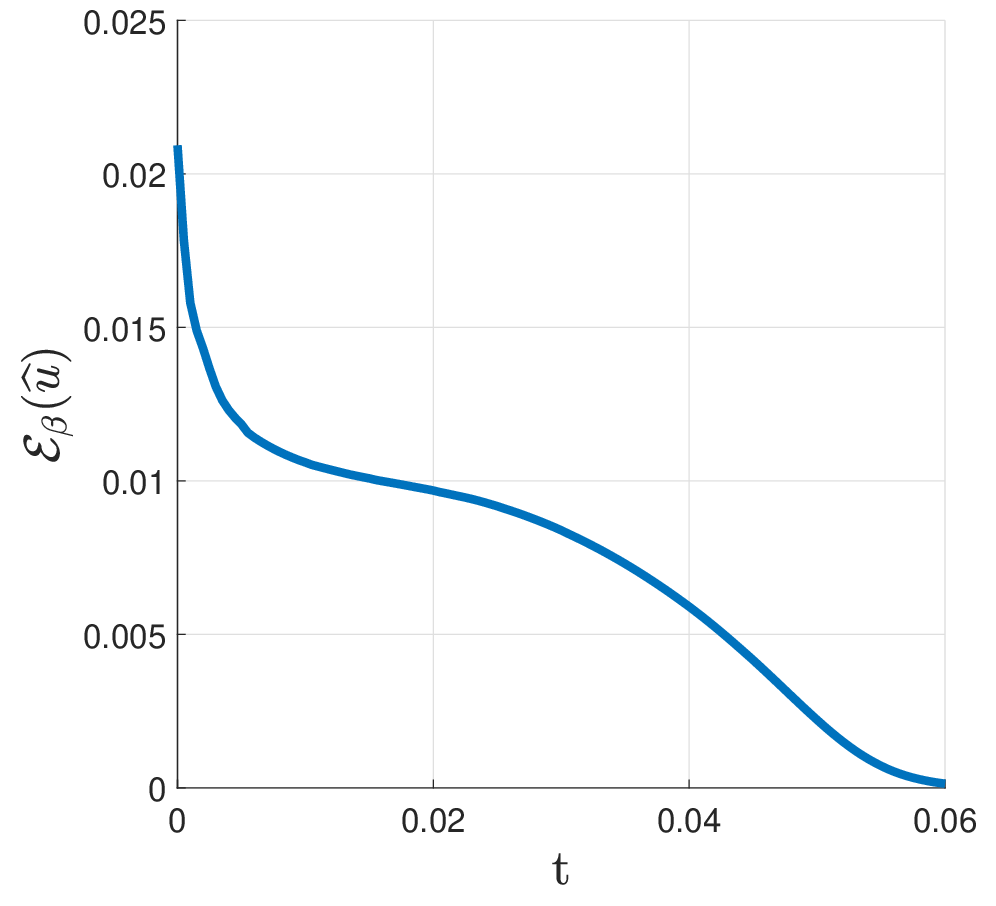}
	\caption{Example~\ref{ex:2D_advection}: Time evolution of the predicted energy functional $\mathcal{E}_{\beta}(\widehat{u})$.}
	\label{fig:2D_advection_energy}
\end{figure}

Last, we consider  an advective Allen--Cahn equation under homogeneous Neumann boundary condition, taken from \cite{RLi_YGao_ZChen_2024}, in the form of 
\begin{eqnarray*}
	u_t &=&  \mu(u) \big( \epsilon^2 \Delta u  - u^3 + u \big) - \beta \cdot \nabla u , \hspace{22mm} \bx  \, \in \, \Omega, \quad t \in [0,T], \\
    \frac{\partial u}{\partial \mathbf{n}} &=& 0,  \hspace{74mm} \bx  \, \in \, \partial \Omega,  \\
	u(\bx, 0) &=& - \tanh \left( \frac{|x_1 + x_2 - 1| + |x_1 - x_2| - 0.1}{\sqrt{2} \epsilon} \right), \hspace{2mm}  \bx \in \Omega,
\end{eqnarray*}
where $\Omega= [0,1]^2$, $\mu(u)=100$, $\epsilon=0.01$, $\beta=(0, -100(x_1 - 0.5))^T$, and $T=0.06$.  Training parameters are given in Table~\ref{tab:prm_ex2D_advection}. We present the evolution of the predicted phase function $\widehat{u}$ in Figure~\ref{fig:2D_advection}. From the results it is seen that the maximum bound principle is satisfied. The similar numerical results were also observed in \cite{RLi_YGao_ZChen_2024}. Due to the advection term, we update the energy function as follows \cite{YWang_XXiao_XFeng_2023}
\begin{equation*}
	\mathcal{E}_{\beta}(u)=\int_{\Omega} \left( \frac{\epsilon^2}{2}|\nabla u|^2 + F(u) + \frac{1}{2}|\beta|^2  \right)dx,
\end{equation*}
where the last term represents the kinetic energy density associated to the vector field $\beta$. Then, the predicted energy function $\mathcal{E}_{\beta}(\widehat{u})$, exhibited in Figure~\ref{fig:2D_advection_energy}, is decreasing over time as expected.

%%%%%%%%%%%%%%%%%%%%%%%%%%%%%%%%%%%%%%%%%%%%%%%%%%%%%%%%%%%%%%%
\section{Conclusion}\label{sec:conclusion}

In this paper, we have proposed a modified physics-informed neural network (PINN) approach to solve the Allen-Cahn equation by imposing the energy dissipation condition as a penalty term into the loss function of the network. A comprehensive numerical experiments including various types of benchmark examples have been conducted to illustrate the outstanding performance of our proposed approach, including its dynamics prediction and generalization ability under different scenarios. To address the challenge of learning random initial conditions, we substitute them with a continuous and smooth analogue using the Fourier series expansion. Numerical results consistently show a decrease in discrete energy, while also highlighting phenomena like phase separation.  Overall, it is evident that deep learning techniques, like PINN, offer appealing alternatives for solving complex engineering challenges.

The code associated with this manuscript is available in a GitHub repository: \\
\href{https://github.com/mustafakutuk/AC\_PINN}{https://github.com/mustafakutuk/AC\_PINN}.

\section*{Acknowledgements}
The authors would like to express their sincere thanks to the anonymous referees for their most valuable suggestions. HY gratefully acknowledges the research support provided by the Scientific and Technological Research Council of Turkey (TUBITAK) under the program TUBITAK 2219-International Postdoctoral Research Fellowship (Project No: 1059B192302207) and would like to thank the Max Planck Institute for Dynamics of Complex Technical Systems, Magdeburg for its excellent hospitality.

\bibliographystyle{unsrt}
%\bibliography{references}  %%% Uncomment this line and comment out the ``thebibliography'' section below to use the external .bib file (using bibtex) .
%\bibliography{C:/Users/hamdullah/Dropbox/Reference/references.bib}
\bibliography{mybib.bib}

\end{document}